\documentclass[12pt]{article}
\usepackage[dvipdfmx]{graphicx}
\usepackage{color}

\usepackage{here}
\usepackage{amsmath,amssymb,amsthm, amscd}
\usepackage{ytableau,mathrsfs}
  
  \makeatletter
  \@addtoreset{equation}{section}
  \makeatother

\renewcommand{\tilde}[1]{\widetilde{#1}}

\newtheorem{dfn}{Definition}[section]
\newtheorem{thm}[dfn]{Theorem}
\newtheorem{lem}[dfn]{Lemma}
\newtheorem{prop}[dfn]{Proposition}

\theoremstyle{remark}
\newtheorem{rem}{Remark}[section]

\begin{document}
\begin{center} 
{\bf {\LARGE Four-dimensional conical symplectic hypersurfaces}}
\end{center} 
\vspace{0.4cm}

\begin{center}
{\large Ryo Yamagishi}
\end{center} 
\vspace{0.4cm}

\begin{abstract}
We show that every indecomposable conical symplectic hypersurface of dimension four is isomorphic to the known one, namely, the Slodowy slice $X_n$ which is transversal to the nilpotent orbit of Jordan type $[2n-2,1,1]$ in the nilpotent cone of $\mathfrak{sp}_{2n}$ for some $n\ge2$. In the appendix written by Yoshinori Namikawa, conical symplectic varieties of dimension two are classified by using contact Fano orbifolds.
\end{abstract}

\section{Introduction}\label{1}
Throughout this paper we work over the complex number field $\mathbb{C}$. A {\em symplectic variety} is a normal variety $X$ such that the regular part $X_\text{reg}$ admits an algebraic symplectic form $\omega$ whose pullback to any resolution of singularities $\tilde{X}\to X$ extends to a regular 2-form on $\tilde{X}$ (cf. \cite{B}). Nilpotent orbit closures and their Slodowy slices in semi-simple Lie algebras are typical examples of symplectic varieties. 

In this paper we only treat {\em conical} ones, that is, affine symplectic varieties with good $\mathbb{C}^*$-actions (see Section \ref{2} for the precise definition). This does not seem to be too restrictive since it is conjectured in \cite[Conjecture 1.8]{K2} that every symplectic variety is (formally) locally a conical one.

It is proved that the moduli space of conical symplectic varieties is discrete \cite{N2}, and thus classification is important. In this paper we study very rare conical symplectic varieties, namely, hypersurfaces. In order to classify conical symplectic hypersurfaces, it is sufficient to classify {\em indecomposable} ones. The followings are the only known indecomposable symplectic hypersurfaces:\\
$\bullet$ (Dimension 2) Kleinian (or ADE) singularities.\\
$\bullet$ (Dimension 4) A series $X_n,\,n\ge2$.\\
$\bullet$ (Dimension 6) A single example $\hat{X}$.\\
The last two were found in \cite{LNS} as Slodowy slices in the nilpotent cones of simple Lie algebras of type $C_n$ and $G_2$ respectively. See Section \ref{2} for the explicit description of $X_n$. In \cite{LNSvS} the authors asked if the three types of the symplectic varieties above are the only indecomposable conical symplectic hypersurfaces. In dimension 2, it is well-known that the only isomorphism types of germs of symplectic varieties are Kleinian singularities (cf. \cite[2.1]{B}), and thus the answer to the question is affirmative in dimension 2. However, we should note that the same Kleinian singularities can admit different structures as conical symplectic varieties. We will give a classification of Kleinian singularities including conical structures in \S3.2 and Appendix. The goal of this paper is to show that the answer to the question is affirmative also in dimension 4. Similarly to the 2-dimensional case, $X_n$ can admit non-unique structures (besides the natural one given in \cite{LNSvS}) as a conical symplectic variety. The main result is the following.

\begin{thm}\label{main}
Let $X$ be an indecomposable conical symplectic hypersurface of dimension four. Then $X$ is isomorphic to $X_n$ as a conical symplectic variety after suitably replacing the $\mathbb{C}^*$-action on $X$.
\end{thm}

The precise definition of an indecomposable conical symplectic variety will be given in Section \ref{2}.

The proof of this theorem is divided into two parts.  Let $X\subset\mathbb{C}^5$ be a 4-dimensional indecomposable conical symplectic hypersurface. Then, each irreducible component $S$ of its singular locus becomes a 2-dimensional conical symplectic variety i.e., a Kleinian singularity or $\mathbb{C}^2$ after taking normalization. In the first part \S4.1 of the proof, we show that some component $S$ itself is in fact normal and has a singularity of type $A$ (Proposition \ref{A}). The main tools to prove this fact are the theory of symplectic hypersurfaces developed in \cite{LNSvS} and the classification of conical symplectic surfaces (\S3.1, 3.3). The procedure to show the existence of a normal type-$A$ surface is roughly explained as follows. The affine space $\mathbb{C}^5=\mathrm{Spec}\,\mathbb{C}[x_1,\dots,x_5]$ admits a natural $\mathbb{C}^*$-action induced from the hypersurface $X\subset\mathbb{C}^5$, and we will see if certain points of the projectivization of the $\mathbb{C}^*$-variety $\mathbb{C}^5$ are in the projectivized singular locus of $X$. If so, we will see that $X$ tends to have a normal type-$A$ surface by using the classification of conical symplectic surfaces (Lemma \ref{p5},\ref{p4} and \ref{p3}). If not so, the Pfaffian condition (see (\ref{Pf})), which was found in \cite{LNSvS}, gives strong restrictions on the possible values of degrees of the generators $x_1,\dots,x_5$. We will also see that $X$ tends to have a normal type-$A$ surface if the degree $d_i$ of some $x_i$ is equal to the weight $s$ of the symplectic form of $X$ (Lemma \ref{*45}). In this way, we will be able to prove the claim with just one exceptional case. Finally, we will exclude this remaining case by direct (computer) calculation.

In the second part \S4.2, we will complete the proof of the main theorem. The idea is to use a coordinate, which will be denoted by $x_{\alpha_1}$, of $\mathbb{C}[x_1,\dots,x_5]$ whose degree is equal to $s$. Such a coordinate exists because of the existence of a normal type-$A$ surface. The symplectic structure of $X$ is translated to a Poisson structure and then extends to the ambient space $\mathbb{C}^5$. This Poisson structure is completely determined by the Poisson matrix $\Theta$. We will determine the possible values of the entries of $\Theta$. For this, we will define a grading $w(-)$ on $\mathbb{C}[x_1,\dots,x_5]$ using $x_{\alpha_1}$, which is an analogue of the weight decomposition of a finite-dimensional $\mathfrak{sl}_2$-module. We will also show that we can make each entry of $\Theta$ homogeneous with respect to $w(-)$ by rechoosing suitable coordinates of $\mathbb{C}[x_1,\dots,x_5]$ (Lemma \ref{wt}). Then, each entry of the Poisson matrix will be doubly homogeneous with respect to $w(-)$ and the original grading, and we can considerably reduce the candidates of the entries. Finally, we will show that $\Theta$ is essentially equal to the Poisson matrix associated to the Slodowy slice $X_n$, which implies that $X$ is isomorphic to $X_n$.

As mentioned above, the classification of 2-dimensional conical symplectic varieties is important for the proof of the main theorem. Yoshinori Namikawa pointed out to the author that the same classification can also be done by using the theory of contact Fano orbifolds, which gives a different perspective for classification of conical symplectic varieties. Namikawa's argument is summarized in a short note written by himself, and he gave permission to the author to use this short note as the appendix to this article.
\\

\noindent
\underline{Convention}\\
In  this paper, whether or not a nonzero monomial $ax_{i_1}^{j_1}\dots x_{i_m}^{j_m}\,(a\in\mathbb{C}^*)$ appears as a term of a polynomial $g\in \mathbb{C}[x_1,x_2,\dots]$ will matter many times.  If the monomial appears in $g$ for some $a\in\mathbb{C}^*$, we simply say that {\em $g$ contains $x_{i_1}^{j_1}\dots x_{i_m}^{j_m}$}.

\vspace{5mm}
\noindent{\bf Acknowledgements}\\
The author would like to thank Yoshinori Namikawa for helpful advice and for giving permission to use his short note for this article. This work was supported by World Premier International Research Center Initiative (WPI), MEXT, Japan, and by JSPS KAKENHI Grant Number JP19K14504.

\section{Conical symplectic varieties as graded Poisson algebras}\label{2}

An affine symplectic variety $X=\mathrm{Spec}\,R$ is called {\em conical} if $R=\bigoplus_{d=0}^\infty R_d$ is a finitely generated graded $\mathbb{C}$-algebra with $R_0=\mathbb{C}$ and the algebraic symplectic form $\omega$ on the regular part $X_\text{reg}$ is homogeneous with respect to the $\mathbb{C}^*$-action $\lambda:\mathbb{C}^*\times X\to X$ coming from the grading on $R$, that is, there exists an integer $s$ such that
$$\lambda(t)^*(\omega)=t^s\omega$$
for any $t\in\mathbb{C}^*$. By \cite[Lemma 2.2]{LNSvS}, the weight $s$ is positive.

Let $X=\mathrm{Spec}\,R$ be a conical symplectic variety. The symplectic form $\omega$ is regarded as an isomorphism $\omega:T_{X_\text{reg}}\to \Omega^1_{X_\text{reg}}$. Then the structure sheaf $\mathcal{O}_{X_\text{reg}}$ admits a natural Poisson structure defined by
$$\{f,g\}=\omega(\omega^{-1}(df),\omega^{-1}(dg))$$
for any sections $f,g$ of $\mathcal{O}_{X_\text{reg}}$. Namely, the pairing
$$\{-,-\}:\,\mathcal{O}_{X_\text{reg}}\times \mathcal{O}_{X_\text{reg}}\to\mathcal{O}_{X_\text{reg}}$$
is a skew $\mathbb{C}$-bilinear form satisfying the Leibniz rule
$$\{fg,h\}=f\{g,h\}+g\{f,h\}$$ 
and the Jacobi identity
$$\{f,\{g,h\}\}+\{g,\{h,f\}\}+\{h,\{g,f\}\}=0$$ 
for any sections $f,g$ and $h$ of $\mathcal{O}_{X_\text{reg}}$.

By the normality of $X$, this uniquely extends to a Poisson structure on $\mathcal{O}_X$. This Poisson bracket is also homogeneous with respect to the $\mathbb{C}^*$-action on $X$, and its weight is $-s$, that is, the Poisson bracket satisfies
$$\{R_k, R_l\}\subset R_{k+l-s}.$$

We assume that $X$ is a hypersurface of dimension $2m$. Let $\bar{x}_1,\dots,\bar{x}_{2m+1}$ be minimal homogeneous generators of $R$ and set $d_i=\text{deg}(\bar{x}_i)$. We may assume that
$$\text{g.c.d.}(d_1,\dots,d_{2m+1})=1$$
by replacing the $\mathbb{C}^*$-action, if necessary. We give a grading on the polynomial ring $P=\mathbb{C}[x_1,\dots,x_{2m+1}]$ by setting $\text{deg}(x_i)=d_i$ so that the surjection $P\to R$ which maps $x_i$ to $\bar{x}_i$ is a graded homomorphism whose kernel is generated by a single (quasi-)homogeneous polynomial $f\in P$. Then by \cite[Lemma 2.7]{LNSvS}, the Poisson structure on $R$ uniquely extends to $P$. The skew-symmetric matrix $\Theta$ whose $(i,j)$-entry $\Theta_{i,j}$ is $\{x_i,x_j\}$ is called the {\em Poisson matrix} of $P$. Note that $\Theta_{i,j}$ is homogeneous of degree $d_i+d_j-s$.

Moreover, the defining equation $f$ of $X$ and $\Theta$ satisfy the following equation after rescaling $f$ \cite[Lemma 2.7]{LNSvS}
\begin{equation}\label{Pf}
\text{grad}(f)=\text{pf}(\Theta).
\end{equation}
In the case $m=2$, this equation is written explicitly as follows.
\begin{equation}\label{4-dim}
\begin{aligned}
\partial f/\partial x_1&=\;\;\;\Theta_{2,3}\Theta_{4,5}-\Theta_{2,4}\Theta_{3,5}+\Theta_{2,5}\Theta_{3,4}\\
\partial f/\partial x_2&=-\Theta_{1,3}\Theta_{4,5}+\Theta_{1,4}\Theta_{3,5}-\Theta_{1,5}\Theta_{3,4}\\
\partial f/\partial x_3&=\;\;\;\Theta_{1,2}\Theta_{4,5}-\Theta_{1,4}\Theta_{2,5}+\Theta_{1,5}\Theta_{2,4}\\
\partial f/\partial x_4&=-\Theta_{1,2}\Theta_{3,5}+\Theta_{1,3}\Theta_{2,5}-\Theta_{1,5}\Theta_{2,3}\\
\partial f/\partial x_5&=\;\;\;\Theta_{1,2}\Theta_{3,4}-\Theta_{1,3}\Theta_{2,4}+\Theta_{1,4}\Theta_{2,3}
\end{aligned}
\end{equation}
These equations will play important roles in proceeding the classification.

\begin{dfn}
A conical symplectic hypersurface $X$ is indecomposable if its unique $\mathbb{C}^*$-fixed point is a Poisson subscheme of $X=\mathrm{Spec}\,R$. Here, a Poisson subscheme means a closed subscheme of $X$ defined by a Poisson ideal $I\subset R$ (i.e. an ideal satisfying $\{I,R\}\subset I$).
\end{dfn}

By \cite[Lemma 2.5]{LNSvS}, in order to classify conical symplectic hypersurfaces, it suffices to classify indecomposable ones.

In Section \ref{4}, we will prove that every indecomposable conical symplectic hypersurface is isomorphic to $X_n$, the only known one. $X_n$ is defined by the following polynomial (cf. \cite{LNSvS}
$$f_n=-ye_0^2+he_0e_1+xe_1^2+\Delta^n$$
of $\mathbb{C}[x,h,y,e_0,e_1]$ where $\Delta=h^2+4xy$. The Poisson matrix is given by
$$\Theta_n=\begin{pmatrix}
0&-2x&h&0&e_0\\
2x&0&-2y&e_0&-e_1\\
-h&2y&0&e_1&0\\
0&-e_0&-e_1&0&2n \Delta^{n-1}\\
-e_0&e_1&0&-2n \Delta^{n-1}&0
\end{pmatrix}.$$
In \cite{LNSvS}, the assigned degrees of $x,h,y,e_0,e_1$ are $2,2,2,2n-1,2n-1$ respectively. However, there are other possibilities of degrees (see (\ref{degree}) in Subsection \ref{4.2}). Thus, even if the $\mathbb{C}^*$-action on $X$ is effective, the replacement of the $\mathbb{C}^*$-action in the statement of Theorem \ref{main} is necessary.

\section{Conical symplectic surfaces}\label{3}

In this section we study conical symplectic surfaces. The results in this section will be used to classify 4-dimensional hypersurfaces in Section \ref{4}.

\subsection{Classification}\label{3.1}
  
In this section we classify 2-dimensional conical symplectic hypersurfaces. It is well-known that the only 2-dimensional symplectic singularities are Kleinian singularities (cf. \cite[2.1]{B}). However, the same singularities can admit different structures as conical symplectic varieties. The complete classification of conical symplectic surfaces including non-singular ones is done by Namikawa in Appendix using the theory of contact Fano orbifolds. As a warm-up for the classification of 4-dimensional hypersurfaces, we give an alternative proof under the assumption that the surfaces are hypersurfaces. Note that this assumption is in fact automatic since Kleinian singularities are hypersurface singularities. In the proof, we will replace homogeneous generators of the coordinate ring by others in order to make the Poisson matrix into a normal form. The similar idea will be used for 4-dimensional hypersurfaces in Section \ref{4}.

Let $X=\{f=0\}$ be a conical symplectic hypersurface in $\mathbb{C}^3$. We take a $\mathbb{C}^*$-equivariant Poisson embedding into $\mathbb{C}^3=\mathrm{Spec}\,\mathbb{C}[x,y,z]$ where $x,y$ and $z$ are homogeneous generators. Let $d_1,d_2$ and $d_3$ be the degrees of $x,y$ and $z$ respectively. Let $s\in\mathbb{N}$ be the weight of the symplectic form $\omega$ on $X_\text{reg}$. By (\ref{Pf}), we have
\begin{equation}\label{deri}
\Theta_{1,2}=\frac{\partial f}{\partial z},\;\Theta_{1,3}=-\frac{\partial f}{\partial y},\;\Theta_{2,3}=\frac{\partial f}{\partial x},
\end{equation}
and in particular the degree of $f\in\mathbb{C}[x,y,z]$ is $d_1+d_2+d_3-s$. One can check that the symplectic form on $X_{\mathrm{reg}}$ corresponding to the Poisson structure $\Theta$ is given by the volume form obtained as the residue
$$Res\left(\frac{dx\wedge dy\wedge dz}{f}\right)=\frac{dx\wedge dy}{\partial f/\partial z}=-\frac{dx\wedge dz}{\partial f/\partial y}=\frac{dy\wedge dz}{\partial f/\partial x}$$
on $X_{\mathrm{reg}}$.

We will prove

\begin{prop}\label{2-dim}(Classification of 2-dimensional conical symplectic hypersurfaces, cf. Proposition in Appendix)
Let $(X,\omega)$ be a 2-dimensional conical symplectic hypersurface with a graded coordinate ring $R$. Assume that the degrees of the homogeneous generators of $R$ are relatively prime. Then $(X,\omega)$ is isomorphic to
$$\left(\mathrm{Spec}\,\mathbb{C}[x,y,z]/(f),\,Res\left(\frac{dx\wedge dy\wedge dz}{f}\right)\right)$$
as a conical symplectic variety where $f$ and the weights $(s,d_1,d_2,d_3)$ of $\omega,x,y$ and $z$ are as in one of the following:\\
$\bullet$ (smooth case) $f=z$, and $d_1$ and $d_2$ are relatively prime integers such that $d_1+d_2=s$.\\
$\bullet$ ($A_n$-type, $n\ge1$) $f=x^{n+1}+yz$, and $d_1,d_2$ and $d_3$ are relatively prime integers satisfying $(n+1)d_1=d_2+d_3$ and $d_1=s$.\\
$\bullet$ ($D_n$-type, $n\ge4$) $f=x^{n-1}+xy^2+z^2$ and $(s,d_1,d_2,d_3)=(1,2,n-2,n-1)$.\\
$\bullet$ ($E_6$-type) $f=x^4+y^3+z^2$ and $(s,d_1,d_2,d_3)=(1,3,4,6)$.\\
$\bullet$ ($E_7$-type) $f=x^3y+y^3+z^2$ and $(s,d_1,d_2,d_3)=(1,4,6,9)$.\\
$\bullet$ ($E_8$-type) $f=x^5+y^3+z^2$ and $(s,d_1,d_2,d_3)=(1,6,10,15)$.
\end{prop}

{\em Proof.} When $X$ is not indecomposable, $X$ falls into the smooth case by \cite[Lemma 2.5]{LNSvS}. 

We assume that $X$ is indecomposable in the rest of the proof. Then $\text{deg}(\Theta_{i,j})$ is nonzero for all $i,j$ (see the proof of Lemma 2.5 of \cite{LNSvS}). Moreover, we have $\Theta_{i,j}\ne0$ for $i\ne j$ \cite[Lemma 2.6 (2)]{LNSvS}. Let $I\subset \mathbb{C}[x,y,z]$ be the ideal generated by $\Theta_{1,2},\Theta_{1,3}$ and $\Theta_{2,3}$. Then $I$ is a nontrivial Poisson ideal and thus the radical $\sqrt{I}$ coincides with the maximal ideal $\mathfrak{m}=(x,y,z)$. (Note that this condition is equivalent to the condition that $X$ has isolated singularity because of (\ref{deri})). In particular, $I$ must contain a power of $z$. This implies that $\Theta_{i,j}$ must contain $z^k$ for some  $i,j$ and $k\in\mathbb{N}$ (see \underline{Convention} in Section \ref{1} for the usage of  the word ``contain" in this paper). We may assume that $d_1\le d_2\le d_3$ by reordering the coordinates. Then $k$ must be 1 since $\deg \Theta_{i,j}=d_i+d_j-s<2d_3$. \\

\noindent
{\bf (Case 1)} $\Theta_{1,2}$ contains $z$.

In this case $\Theta_{1,2}$ is of the form $az+p(x,y)$ for some $a\in\mathbb{C}^*$ and $p(x,y)\in\mathbb{C}[x,y]$ since $d_i$'s are positive. Thus, we may assume that $\Theta_{1,2}=2z$ by using $\frac{1}{2}(az+p(x,y))$ as a new coordinate instead of $z$. Similarly, $\Theta_{1,3}$ or $\Theta_{2,3}$ must contain $y^l$ for some $l\in\mathbb{N}$. Then $l=1$ or $2$ since we have
$$\deg \Theta_{i,j}=d_i+d_j-s<d_2+d_3=d_1+2d_2-s<3d_2$$
noticing that $\deg(\Theta_{1,2})=\deg(z)$.\\

\noindent
{\bf (Case 1-1)} $\Theta_{1,3}$ or $\Theta_{2,3}$ contains $y^2$.

Since $d_3=d_1+d_2-s$ and $d_1\le d_2$, we have
$$\deg(z^2)=2d_3=2d_1+2d_2-2s<2d_1+2d_2=\deg(x^2 y^2)\le d_1+3d_2=\deg(x y^3).$$
Combining this with (\ref{deri}), we can write
$$f=z^2+a_1 y^3+a_2 xy^2+a_3 x^m y+a_4x^n$$
with some $a_i\in\mathbb{C}$ such that $a_1\ne0$ or $a_2\ne0$ and $m,n\in\mathbb{N}$.\\

\noindent
{\bf (Case 1-1-1)} $a_1\ne0$.

We may assume that $a_2=0$ by using a new coordinate $y+\frac{a_2}{3a_1}x$ instead of $y$ if $a_2\ne0$, which happens only when $d_1=d_2$ since $f$ is homogeneous. Then we have
$$(d_1,d_2,d_3)=(d_1,2d_1-2s,3d_1-3s)$$
and in particular $d_1-2s=d_2-d_1\ge0$. We also have $(a_3,a_4)\ne(0,0)$ since otherwise $\partial f/\partial x=\Theta_{2,3}=0$. We first consider the case $a_3\ne0$. Then we see that $m=2$ or $3$ since
$$m d_1=\deg(x^m)=\deg(y^2)=4d_1-4s$$
and $d_1\ge 2s$. If $m=3$, then we have $(d_1,d_2,d_3)=(4s,6s,9s)$ and we can write
$$f=z^2+a_1 y^3+a_3 x^3 y$$
with $a_1 a_3\ne0$ since $d_1\nmid \deg(f)$. This $f$ is the equation for $E_7$ after rescaling $x$ and $y$. If $m=2$, then $(d_1,d_2,d_3)=(2s,2s,3s)$ and $q(x,y)=f-z^2\in\mathbb{C}[x,y]$ is a homogeneous cubic polynomial. The condition that $X$ has isolated singularity implies that the equation $q(x,y)=0$ has mutually distinct 3 solutions $[x\mathbin:y]\in\mathbb{P}^1$, and therefore $q(x,y)$ becomes $x^3+xy^2$ after a suitable linear coordinate change. Then $f$ is the equation for $D_4$.

Next we consider the case $a_4\ne0$. Then we see that $n=3,4$ or $5$ since
$$n d_1=\deg(x^n)=\deg(y^3)=6d_1-6s$$
and $d_1\ge 2s$. If $n=3$, then $f$ can be transformed into the equation for $D_4$ similarly to the case $a_3\ne0$. If $n=4$, then  we have $(d_1,d_2,d_3)=(3s,4s,6s)$ and we can write
$$f=z^2+a_1 y^3+a_4 x^4$$
with $a_1 a_4\ne0$. This $f$ is the equation for $E_6$ after rescaling $x$ and $y$. Finally, if $n=5$, then  we have $(d_1,d_2,d_3)=(6s,10s,15s)$ and we can write
$$f=z^2+a_1 y^3+a_4 x^5$$
with $a_1 a_4\ne0$. This $f$ is the equation for $E_8$ after rescaling $x$ and $y$.\\

\noindent
{\bf (Case 1-1-2)} $a_2\ne0$.

We may assume that $a_1=0$ since the case $a_1\ne0$ has just been treated above. Then we may also assume that $a_3=0$ by using $y+\frac{a_3}{2a_2}x$ as a new coordinate instead of $y$. Then we have $a_4\ne0$ since otherwise any of the partial derivatives of $f$ would not contain a power of $x$ and hence $\sqrt{I}\ne\mathfrak{m}$. Therefore, we can write
$$f=z^2+a_2 xy^2+a_4x^n$$
with $a_2, a_4\ne0$. This $f$ is the equation for $D_{n+1}$ after rescaling $x$ and $y$.

Since $\deg(z^2)=\deg(xy^2)=\deg(x^n)$, we have $2(d_1+d_2-s)=d_1+2d_2=n d_1$ and thus obtain
$$(d_1,d_2,d_3)=(2s,(n-1)s,ns)$$
with $n\ge3$.\\

\noindent
{\bf (Case 1-2)} $\Theta_{1,3}$ or $\Theta_{2,3}$ contains $y$.

If $\Theta_{1,3}$ contains $y$, then we can write $\Theta_{1,3}=-\frac{\partial f}{\partial y}=ay+bx^m$ for some $a\in\mathbb{C}^*,b\in\mathbb{C}$ and $m\ge1$ since $\frac{\partial f}{\partial z}=2z$. We may assume that $\Theta_{1,3}=ay$ by using a new coordinate $y+\frac{b}{a}x^m$ instead of $y$. Note that this replacement does not break the condition $\Theta_{1,2}=2z$. Then we can write
$$f=cx^{n+1}-\frac{a}{2}y^2+z^2$$
for some $c\in\mathbb{C}$ and $n\ge1$. We see that $c$ is nonzero by the condition $\sqrt{I}=\mathfrak{m}$, and thus $X$ isomorphic to $A_n$-type with $d_2=d_3$. If $\Theta_{2,3}$ contains $y$, then we have $\deg(\Theta_{2,3})=d_2+d_3-s=d_2$ and hence $d_3=s$. In this case we have $d_1=d_2=s$ since
$$s=\deg(z)=\deg(\Theta_{1,2})=d_1+d_2-s$$
and $d_1\le d_2\le d_3=s$. Therefore, we can write $f=z^2+a_1y^2+a_2xy+a_3 x^2$. The condition that $X$ has isolated singularity implies that the quadratic equation $a_1y^2+a_2xy+a_3 x^2=0$ has mutually distinct 2 solutions $[x\mathbin:y]\in\mathbb{P}^1$, and therefore $a_1y^2+a_2xy+a_3 x^2$ becomes $xy$ after a suitable linear coordinate change. Then $f$ is the equation for $A_1$. \\

\noindent
{\bf (Case 2)} $\Theta_{1,2}$ does not contain $z$.

In this case $\Theta_{1,2}$ is of the form $ax^k+by$ for some $k\in\mathbb{N}$ and $a,b\in\mathbb{C}$ for a degree reason. If $b\ne0$, we may assume that $\Theta_{1,2}=y$ by using new coordinates $\frac{x}{b}$ and  $y+\frac{a}{b}x^k$ instead of $x$ and $y$ respectively. Then $\Theta_{1,3}$ is of the form
$$-z+\sum_{k,l} a_{k,l}x^k y^l$$
by (\ref{deri}). We may assume that $\Theta_{1,3}=-z$ by using a new coordinate $z-\sum_{k,l} \frac{a_{k,l}}{l+1} x^k y^l$ instead of $z$. Indeed, we have
$$\begin{aligned}
\left\{x,z-\sum_{k,l} \frac{a_{k,l}}{l+1} x^k y^l\right\}&=\{x,z\}-\left\{x,\sum_{k,l} \frac{a_{k,l}}{l+1} x^k y^l\right\}\\
&=-z+\sum_{k,l} a_{k,l}x^k y^l-\sum_{k,l} \frac{la_{k,l}}{l+1} x^k y^l\\
&=-\left(z-\sum_{k,l} \frac{a_{k,l}}{l+1} x^k y^l\right).
\end{aligned}$$
Then we can write
$$f=cx^{n+1}+yz$$
for some $c\in\mathbb{C}^*$, and $X$ isomorphic to $A_n$-type. If $b=0$, then we have
$$f=ax^kz+p(x,y)$$
for some $p(x,y)\in\mathbb{C}[x,y]$. One sees that $k=1$ from the condition that $z$ is contained in $\sqrt{I}$. This case can be reduced to the previous case $\Theta_{1,2}=y$ by switching $x$ and $y$.
\qed

\vspace{3mm}

\begin{rem}
In the proof above, the condition $s>0$ plays crucial roles. Indeed, this condition is essential for the classification of conical symplectic surfaces since there are lots of normal surface singularities other than ADE singularities whose coordinate rings admit graded Poisson structures. For example, every isolated surface singularity defined by a quasi-homogeneous polynomial admits a graded Poisson structure on the coordinate ring which is induced from the volume form on the smooth part of the surface.
\end{rem}

\subsection{Non-normal surfaces}\label{3.2}

In this subsection we study surfaces whose normalizations are conical symplectic surfaces. This will be necessary later since such surfaces could appear in 4-dimensional symplectic hypersurfaces (or more generally complete intersections). In general, symplectic varieties admit stratification by locally closed connected symplectic submanifolds, and the closure of each stratum becomes a symplectic variety after taking normalization (cf. \cite[Theorem 2.4]{K1}). When a symplectic variety is a complete intersection, the singular locus is of codimension at most 3 \cite[Proposition 1.4]{B}. Thus, in particular, 4-dimensional symplectic hypersurface contains a subvariety whose normalization is smooth or has an ADE-singularity. When the 4-dimensional hypersurface is conical, each irreducile component $S$ of the singular locus is stable under the $\mathbb{C}^*$-action. This action extends to the normalization $\tilde{S}$ of $S$ (see e.g. \cite[Theorem 2.10]{P}). The Poisson structure of $S$ also extends to $\tilde{S}$ (cf. \cite{K3}), and $\tilde{S}$ becomes a conical sympectic surface. Note that the coordinate ring $A=H^0(\mathcal{O}_S)$ is positively graded and so is $B=H^0(\mathcal{O}_{\tilde{S}})$. Note also that the normalization map $\tilde{S}\to S$ is a bijection since both $\tilde{S}$ and $S$ have the unique $\mathbb{C}^*$-fixed points.

Possible structures of $B$ are classified in the previous subsection. Note that, if $B$ is smooth or of type $A_n$, then $B$ is isomorphic to a graded Poisson subalgebra of the polynomial ring $\mathbb{C}[u,v]$ of two variables where the Poisson structure is given by $\{u,v\}=1$. In the case of type $A_n$, $B$ is isomorphic to $\mathbb{C}[u^{n+1},v^{n+1},uv]$ as a graded Poisson algebra by setting
$$\deg(u)=\frac{d_1}{n+1}\;\;\text{   and   }\;\;\deg(v)=\frac{d_2}{n+1}.$$
The smooth case is regarded as $A_0$-type.

Let $-s$ be the weight of the Poisson bracket of $A$ (and hence of $B$). The following lemma shows that $A$ is isomorphic to a monomial algebra if $A$ contains a nonzero element of degree $s$.

\begin{lem}\label{mono}
Let $A$ be a graded Poisson algebra as above. If $A$ contains a nonzero homogeneous element of degree $s$, then the integral closure $B$ is of type $A_n$ ($n\ge0$). Moreover, $A$ is isomorphic to a graded Poisson subalgebra generated by monomials of $\mathbb{C}[u,v]$ where the grading on $\mathbb{C}[u,v]$ is the one defined as above.
\end{lem}

{\em Proof.} The first claim is clear since, for types $D_n$ and $E_n$ in Proposition \ref{2-dim}, the algebras of functions (and hence their subalgebras) are generated by elements whose degrees are greater than $s$.

For the second claim, we regard $A$ as a subalgebra of $\mathbb{C}[u,v]$. We may assume that $\deg(u)\le\deg(v)$. Let $h\in A\subset \mathbb{C}[u,v]$ be a nonzero homogeneous element of degree $s$. We first show that we may assume that $h$ is a monomial by replacing homogeneous generators of $\mathbb{C}[u,v]$. If $\deg(u)<\deg(v)$, then $h$ is of the form $auv+bu^k$ for some $a,b\in\mathbb{C}$ and $k>1$. If $a=0$, we are done. If $a\ne0$, then we can make $h$ into a monomial by using a new coordinate $v-\frac{b}{a}u^{k-1}$ instead of  $v$. Note that this replacement preserves the Poisson bracket $\{u,v\}=1$. If $\deg(u)=\deg(v)$, then $h$ is of the form
$$au^2+buv+cv^2$$
for some $a,b,c\in\mathbb{C}$. We can make $h$ into a monomial by a linear transformation of $SL(2)$ since $h$ is a quadratic form of rank at most 2. Note again that this operation preserves the Poisson structure.

Now we can assume that $h$ is $uv$ or $u^k,\,k\ge2$. First we assume $h=uv$. Let $p\in A$ be any nonzero homogeneous element and let
$$p=p_1+\cdots+p_m$$
be the decomposition into monomials $p_i=c_i u^{a_i}v^{b_i}$ with $c_i\in\mathbb{C}^*$. We will show that each $p_i$ is in $A$. We may assume that every $a_i-b_i$ is nonzero since $uv\in A$. Then $a_i-b_i$ are different for all $i$ since
$$\deg(p)=\deg (u)a_i+\deg (v)b_i=\deg (u)(a_i-b_i)+(\deg (u)+\deg (v))b_i$$
is constant by the homogeneity of $p$. We define $p^{(j)},\,j=1,\dots,m$ inductively by setting $p^{(1)}=p$ and
$$p^{(j+1)}=\{p^{(j)},uv\}\in A.$$
Note that $\{u^{a_i}v^{b_i},uv\}=(a_i-b_i)u^{a_i}v^{b_i}$ and therefore
$$p^{(j)}=\sum_{i=1}^m (a_i-b_i)^{j-1} c_i p_i.$$
We can show that $p^{(1)},\dots,p^{(m)}$ are linearly independent over $\mathbb{C}$. Indeed, if
$$\sum_{j=1}^m e_j p^{(j)}=0$$
for $e_j\in\mathbb{C}$, then we have
$$\sum_{i=1}^m \sum_{j=1}^m e_j(a_i-b_i)^{j-1} c_i p_i=0.$$
The linear independence of $c_i p_i$'s implies that
$a_1-b_1,\dots,a_m-b_m$ are $m$ distinct solutions of the equation
$$e_1+e_2 t+\cdots+e_m t^{m-1}=0$$
and hence all $e_j$ must be zero. This implies that each $p_i$ is obtained as a $\mathbb{C}$-linear combination of $p^{(j)}$'s and is in particular an element of $A$.

When $h=u^k$, we have $\deg(v)=s-\deg(u)=(k-1)\deg(u)$. Thus, any nonzero homogeneous element of $A$ is written as
$$p=\sum_{i=1}^m c_i u^{l+(k-1)i} v^{m-i}$$
for some $l,m\ge0$ and $c_i\in\mathbb{C}$ with $c_1\ne0$. Similarly as above, we define $p^{(j)},\,j=1,\dots,m$ inductively by $p^{(1)}=p$ and $p^{(j+1)}=\{p^{(j)},u^k\}\in A$. Then $p^{(j)}$'s are linearly independent since $p^{(j)}$ consists of terms $au^b v^c$ with $c\le m-j$ and $a\ne0$ when $c=m-j$. This shows that $p^{(j)}$'s span the same vector space as the $m$-dimensional space spanned by $u^{l+(k-1)i} v^{m-i}$'s and in particular $u^{l+(k-1)i} v^{m-i}\in A$.
\qed

\vspace{3mm}

The following lemma will be very useful in proceeding classification of 4-dimensional symplectic hypersurfaces.


\begin{lem}\label{non-normal}
Let $A$ be a graded Poisson algebra which is the coordinate ring of an irreducible component $S$ of a 4-dimensional conical symplectic hypersurface $(X,\omega)$. Assume that $A$ contains a nonzero homogeneous element of the same degree $s>0$ as the weight of the symplectic form $\omega$. Then $A$ is integrally closed.
\end{lem}

{\em Proof.} By Lemma \ref{mono}, we may assume that the integral closure $B$ of $A$ is a Poisson subalgebra of $\mathbb{C}[u,v]$ generated by $x=u^{n+1},y=v^{n+1}$ and $z=uv$ possibly with $n=0$. We may also assume that $A$ is generated by monomials of $\mathbb{C}[u,v]$. Since the normalization $\mathrm{Spec}(B)\to \mathrm{Spec}(A)$ is a bijection, we have $\sqrt{\mathfrak{m}B}=(x,y,z)$ in $B$ where $\mathfrak{m}$ is the maximal ideal of $A$ corresponding to the unique $\mathbb{C}^*$-fixed point. In particular, the set of generators of $A$ must contain powers of $x$ and $y$.

Let $k$ and $l$ be the smallest positive integers such that $x^k,y^l\in A$. For any monomial $u^a v^b\in A$, we have $n+1|a-b$ since $A\subset B$. In fact we can choose $u^a v^b\in A$ such that $a-b=n+1$ since otherwise $x=u^{n+1}$ would not be contained in the integral closure $B$. By successively applying $\{x^k,-\}$ to this monomial, we obtain $x^{k'}\in A$ such that $k$ and $k'$ are coprime. Indeed, we have
$$\{x^k,u^a v^b\}=k(n+1)b u^{k(n+1)+a-1}v^{b-1},$$
and $k$ and $\frac{k(n+1)+a-1-(b-1)}{n+1}=k+\frac{a-b}{n+1}$ are coprime if $k$ and $\frac{a-b}{n+1}$ are coprime. Similarly, we obtain $y^{l'}\in A$ such that $l$ and $l'$ are coprime. We choose the smallest $k'$ (resp. $l'$) such that $x^{k'}\in A$ (resp. $y^{l'}\in A$) and that $k$ and $k'$ (resp. $l$ and $l'$) are coprime.

We consider the two cases (1) $A$ contains $z=uv$, and (2) $A$ does not contain $z=uv$, separately. First we assume (1) $z=uv\in A$.

If $A$ contains $x$ and $y$, then $A=B$ and hence we are done. If $k,l\ge2$, then $A$ must contain $x^k,x^{k'},y^l,y^{l'}$ and $z$ as part of minimal monomial generators but they must generate whole $A$ since $A$ is generated by at most 5 elements. However, the algebra generated by these 5 elements is not closed under the Poisson bracket. Indeed, one can check that $A$ does not contain $\{x^k,y^l\}$ if $k\ne l$ and otherwise $A$ does not contain $\{x^{k'},y^l\}$ noticing that $k$ and $k'$ are coprime.

We may assume that $k\ge2$ and $l=1$ by switching $u$ and $v$ if necessary. In this case $A$ must contain $x^k,x^{k'},y$ and $z$ as part of minimal monomial generators. By applying $\{-,y\}$ successively to $x^k$, we see that $A$ must contain $xz^m$ for some $m\ge1$ in order for $A$ to be closed under the Poisson bracket. We choose the smallest $m$ such that $xz^m\in A$. For $\{x^k, xz^m\}$ to be in $A$, we must have $k'=k+1$. Similarly, we must have $k=2$ for $\{x^{k+1}, xz^m\}$ to be in $A$. When $k=2$ and $k'=3$, the algebra generated by
$$x^2,x^3,y,z\text{ and }xz^m$$
is indeed a Poisson algebra when $m\le n$. Since the Poisson structure on $A=\mathbb{C}[z,xz^m,x^2,x^3,y]$ comes from a hypersurface $\{f=0\}\subset\mathbb{C}^5$, it must lift to $P:=\mathbb{C}[x_1,x_2,x_3,x_4,x_5]$ as explained in Section \ref{2}. We show that this is impossible.

Let $\phi:P\to A$ be the graded Poisson surjection defined by
$$(x_1,\dots,x_5)\mapsto (z,xz^m,x^2,x^3,y).$$
We define a new $\mathbb{Z}$-grading $w(-)$ on $P$ by setting
$$w(x_1)=0,\,w(x_2)=-(n+1),\,w(x_3)=-2(n+1),\,w(x_4)=-3(n+1),\,w(x_5)=n+1$$
so that $\{z,\phi(x_i)\}=w(x_i)\phi(x_i)$ for every $i$. (See Introduction or Subsection \ref{4.2} for the motivation to introduce the grading $w(-)$.) We say that a polynomial $g\in P$ is $w$-{\em homogeneous} if $g$ is homogeneous with  respect to $w(-)$. Then we may assume that $\{x_i,x_j\}$ is $w$-homogeneous of degree $w(x_i)+w(x_j)$ for any $i,j$ by suitably replacing the homogeneous generators $x_1,\dots,x_5$ of $P$ with others. This follows from the same argument as in the proof of Lemma \ref{wt} in Subsection \ref{4.2}. The equation (\ref{Pf}) shows that $\{x_1,f\}=0$, and this implies that $f$ is $w$-homogeneous and $w(f)=0$. However, by (\ref{4-dim}), we have
$$w(f)=\sum_{i=1}^5 w(x_i)=-5(n+1)\ne0,$$
which is a contradiction.

Next we assume (2) $z=uv\not\in A$. We show that this case is impossible. Let $\bar{x}_1,\dots,\bar{x}_5\in A\subset\mathbb{C}[u,v]$ be monic monomials which generate $A$, and let $I\subset A$ be the ideal generated by all the Poisson brackets $\{\bar{x}_i,\bar{x}_j\}$. Since $I$ is clearly a nontrivial Poisson ideal and $\mathrm{Spec}\,A$ has isolated singularity, the radical $\sqrt{I}$ is equal to the maximal ideal $\mathfrak{m}$ of $A$. In particular, $I$ contains a power of $x$. For this, some $\bar{x}_i$ must be equal to $x^{k''} z$ for some $k''\ge0$. Indeed, in order for
$$\{u^a v^b,u^{a'}v^{b'}\}=(a b'-a' b)u^{a+a'-1}v^{b+b'-1}$$
to be a nonzero multiple of a power of $x=u^{n+1}$, one of $b$ and $b'$ must be equal to $1$. Since $z=uv\not\in A$ by assumption, we have $k''\ge1$. Similarly, some $\bar{x}_j$ must be equal to $y^{l''} z$ with $l''\ge1$. We choose the smallest $k''$ and $l''$ satisfying $x^{k''} z,y^{l''} z\in A$. Then the monomials
$$x^k,y^l,x^{k''} z,y^{l''} z$$
must be part of minimal monomial generators of $A$. ($x^{k'}$ (resp. $y^{l'}$) is also necessary as a generator if $k\ge2$ (resp. $l\ge2$).)  In this case we see that $A$ needs more than 5 generators in order for $A$ to be closed under the Poisson brackets, which is a contradiction.
\qed

\subsection{Ramification points on projectivizations}

Conical symplectic surfaces are $\mathbb{C}^*$-varieties and thus we can consider their projectivizations. In this subsection we review such projectivizations.

Let $S=\mathrm{Spec}\,A$ be a conical symplectic surface. The projectivization $\mathbb{P}(S)$ is nothing but $\mathrm{Proj}\, A$. As a set, $\mathbb{P}(S)$ consists of the $\mathbb{C}^*$-orbits in $S$ minus the origin. Note that, if the g.c.d. of the degrees of nonzero generators of $A$ is $d$, then the stabilizer subgroup of $\mathbb{C}^*$ for a general point of $S$ has order $d$. We call $p\in \mathbb{P}(S)$ a {\em ramification point} if the order $\mathrm{Stab}_p$ of the stabilizer group of some (and hence any) point of the $\mathbb{C}^*$-orbit in $S$ corresponding to $p$ is greater than $d$. We define the {\em ramification index} $r_p$ of $p$ as $\mathrm{Stab}_p/d$, which is always a positive integer.

\vspace{3mm}

\begin{rem}
In general projectivizations of conical symplectic varieties have contact Fano orbifold structures, and conversely we can recover conical symplectic varieties from contact Fano orbifolds. For details, see Appendix and the references therein. In the surface case $\mathbb{P}(S)$ is isomorphic to $\mathbb{P}^1$ as a scheme but there are various orbifold structures on $\mathbb{P}^1$. In Appendix, the classification of conical symplectic surfaces is done by classifying possible Fano orbifold structures.
\end{rem}

From the concrete description of conical symplectic surfaces in Proposition \ref{2-dim}, we can find all the ramification points of $\mathbb{P}(S)$.

\begin{prop}(cf. Lemma in Appendix)\label{ram}
Let  $S=\mathrm{Spec}\,\mathbb{C}[x,y,z]/(f)$ be a conical symplectic surface with $d_1:=\deg(x),d_2:=\deg(y)$ and $d_3:=\deg(z)$. Then the ramification points of $\mathbb{P}(S)$ are listed as follows.\\
$\bullet$ (smooth case) $f=z$; the points $[0\mathbin:1\mathbin:0]$ and $[1\mathbin:0\mathbin:0]$ with ramification index $d_1/d$ and $d_2/d$ respectively as long as $d_1,d_2>d:=\mathrm{g.c.d.}(d_1,d_2)$.\\
$\bullet$ ($A_n$-type, $n\ge1$) $f=x^{n+1}+yz$; the points $[0\mathbin:1\mathbin:0]$ and $[0\mathbin:0\mathbin:1]$ with ramification index $d_2/d$ and $d_3/d$ respectively as long as $d_2,d_3>d:=\mathrm{g.c.d.}(d_1,d_2,d_3)$.\\
$\bullet$ ($D_n$-type, $n\ge4$) $f=x^{n-1}+xy^2+z^2$; the point $[0\mathbin:1\mathbin:0]$ with ramification index $n-2$, and the two points $[1\mathbin:\pm\zeta_4\mathbin:0]$ (resp. $[1\mathbin:0\mathbin:\pm\zeta_4]$) if $n$ is even (resp. odd) with ramification index $2$.\\
$\bullet$ ($E_6$-type) $f=x^4+y^3+z^2$; the three points $[0\mathbin:1\mathbin:\zeta_4]$, $[1\mathbin:0\mathbin:\zeta_4]$ and $[1\mathbin:0\mathbin:-\zeta_4]$ with ramification index $2,3$ and $3$ respectively.\\
$\bullet$ ($E_7$-type) $f=x^3y+y^3+z^2$; the three points $[1\mathbin:-1\mathbin:0]$, $[0\mathbin:1\mathbin:\zeta_4]$ and $[1\mathbin:0\mathbin:0]$ with ramification index $2,3$ and $4$ respectively.\\
$\bullet$ ($E_8$-type) $f=x^5+y^3+z^2$; the three points $[1\mathbin:\zeta_6\mathbin:0]$, $[1\mathbin:0\mathbin:\zeta_4]$ and $[0\mathbin:1\mathbin:\zeta_4]$ with ramification index $2,3$ and $5$ respectively.\\
Here $\zeta_i$ is the $i$-th primitive root of unity and $[a\mathbin:b\mathbin:c]$ is the point of the weighted projective space $\mathbb{P}(d_1,d_2,d_3)$ corresponding to a point $(a,b,c)\in S\subset\mathbb{C}^3=\mathrm{Spec}\,\mathbb{C}[x,y,z]$.
\end{prop}

\vspace{3mm}

\begin{rem}
In Proposition \ref{ram}, the $\mathbb{C}^*$-action on $S$ need not be effective, or in other words the degrees of the generators of the coordinate rings need not be relatively prime. The order $\mathrm{Stab}_p$ of the stabilizer group of a ramification point $p\in S$ depends on the g.c.d. of the degrees, but the ramification index $r_p$ does not.
\end{rem}

\section{4-dimensional hypersurfaces}\label{4}

In this section we give a proof of Theorem \ref{main}. The proof is divided into two parts. The first part (Subsection \ref{4.1}) is devoted to show that the 4-dimensional hypersurface $X$ has a normal surface of type $A_k$ as an irreducible component of the singular locus of $X$. In the second part (Subsection \ref{4.2}), we complete the proof of Theorem \ref{main}.

\subsection{$X$ contains a conical symplectic surface of type $A_k$}\label{4.1}

Let $X=\{f=0\}\subset \mathrm{Spec}\,R$ be the hypersurface with $R=\mathbb{C}[x_1,x_2,x_3,x_4,x_5]$ as in section \ref{2}. In particular, we assume that $d_1\le\dots\le d_5$ and $\text{g.c.d.}(d_1,\dots,d_5)=1$. Note that the order of the stabilizer group of a point $p\in X$ with respect to the $\mathbb{C}^*$-action on $X$ is given by
$$\mathrm{Stab}_p=\mathrm{g.c.d.}\{d_i\mid p\not\in \{x_i=0\}\subset\mathbb{C}^5\}.$$

We also assume that $X$ is indecomposable. This is equivalent to the condition that every irreducible component $S$ of the (reduced) singular locus of $X$ is singular. Indeed, if $X$ is indecomposable, then the $\mathbb{C}^*$-fixed point of $X$ is a nontrivial Poisson subscheme of $S$, which implies that $S$ is singular \cite[Lemma 1.4]{K1}.

This subsection is devoted to prove the following.

\begin{prop}\label{A}
$X$ contains a conical symplectic surface of type $A_k\,(k\ge1)$ as an irreducible component of the reduced singular locus of  $X$.
\end{prop}

Let $\mathbb{P}:=\mathbb{P}(d_1,\dots,d_5)=\mathrm{Proj}\,\mathbb{C}[x_1,x_2,x_3,x_4,x_5]$ be the weighted projective space. Note that the reduced singular locus $S'\subset X$ is stable under the $\mathbb{C}^*$-action and we can consider the projectivizations $\mathbb{P}(S')\subset \mathbb{P}(X)\subset\mathbb{P}$. We consider the points
$$p_i=[\delta_{i,1}\mathbin:\delta_{i,2}\mathbin:\delta_{i,3}\mathbin:\delta_{i,4}\mathbin:\delta_{i,5}]\in\mathbb{P}$$
for $i=1,\dots,5$ where $\delta_{i,j}$ is the Kronecker delta.

\begin{lem}\label{p5}
If $p_5$ is in $\mathbb{P}(S')$, then $S'$ contains an irreducible component $S$ which is a normal symplectic surface of type $A_k$ for some $k\ge1$.
\end{lem}

{\em Proof.} Let $S$ be the irreducible component of $S'$ such that $p_5\in\mathbb{P}(S)$ and let $I_S\subset R$ be the defining ideal of $S$. Note that such $S$ is unique since irreducible components of $S'$ intersect only at the origin. For an element $g\in R$, we denote the image of $g$ in $R/I_S$ by $\bar{g}$. Note that the normalization $\tilde{S}$ of $S$ is a conical symplectic surface and that $\tilde{S}$ and $S$ are isomorphic as $\mathbb{C}^*$-varieties outside the singular points (and hence $\mathbb{P}(S)\cong \mathbb{P}(\tilde{S})$ as orbifolds). Therefore, the corresponding $\mathbb{C}^*$-orbits of $\tilde{S}$ and $S$ have the same stabilizer groups, and, in particular, ramification points of $\mathbb{P}(S)$ bijectively correspond to those of $\mathbb{P}(\tilde{S})$ preserving ramification indices.

We first consider the case where $p_5\in\mathbb{P}(S)$ is not a ramification point. Note that the coordinate ring $R/I_S$ of $S$ is generated by at least 3 elements since $S$ is singular by the assumption that $X$ is indecomposable. We first assume that $R/I_S$ is generated by 3 homogeneous elements, then $S$ is normal by Serre's normality criterion. Then we can choose 3 elements from $\bar{x}_1,\dots,\bar{x}_5$ as generators of $R/I_S$ since $R/I_S$ is a positively graded ring. Since $p_5\in\mathbb{P}(S)$, we can choose $\bar{x}_5$ as one of 3 minimal generators of $R/I_S$. Then, from Proposition \ref{2-dim} and the maximality of $d_5$, we see that only $A_k$-type is possible in order for the condition $p_5\in\mathbb{P}(S)$ to be satisfied. Indeed, for type $D_k$ and $E_k$, the vanishing of the two generators $x$ and $y$ of lower degrees in the defining equation of $S$ gives the vanishing of $z$.

If $R/I_S$ is generated by at least 4 elements, general points of $S$ are contained in at most one coordinate hyperplane $\{x_i=0\}\subset \mathbb{C}^5$. Since $p_5$ is assumed not to be a ramification point, general points of $S$ must have stabilizer groups of order $d_5$ with respect to the $\mathbb{C}^*$-action on $\mathbb{C}^5$. Thus, at least 3 of 4 integers $d_1,d_2,d_3$ and $d_4$ are multiples of $d_5$. By the maximality of $d_5$, we have
$$d_2=d_3=d_4=d_5=:\alpha.$$
If $\alpha\ne s$, then we have
$$\alpha\ne\deg{\Theta_{i,j}}=2\alpha-s<2\alpha$$
for $i,j\ge2$ and in particular $\Theta_{i,j}\in \mathbb{C}[x_1]$ for $i,j\ge2$. This shows $x_1|f$ by the Pfaffian condition (\ref{Pf}), which is contrary to the irreducibility of $f$. Thus we obtain $\alpha=s$. If $d_1<\alpha=s$, then we would have $\Theta_{1,j}\in \mathbb{C}[x_1]$ for $j\ge1$, which leads to the same contradiction. Now we have that every $d_i$ is equal to $s$, and thus we have $s=1$. Then, $X$ is homogeneous in the sense of \cite{N1}, and the main result of \cite{N1} states that the only homogeneous singular symplectic varieties of complete intersection are nilpotent cones in semisimple Lie algebras. The only 4-dimensional nilpotent cone is the product of two copies of $A_1$-singularity, but this is not a hypersurface.

Next we assume that $p_5\in\mathbb{P}(S)$ is a ramification point with ramification index $r_{p_5}>1$. We show that $\tilde{S}$ must be smooth or of type $A_k$ using the classification (Proposition \ref{2-dim} and \ref{ram}). Note that $R/I_S$ is generated by homogeneous elements whose degrees are less than or equal to $d_5=d r_{p_5}$ where $d$ is the order of the stabilizer group of general points of $S$. If the integral closure $H^0(\mathcal{O}_{\tilde{S}})$ of $R/I_S$ were isomorphic to the coordinate ring $R_{E_k}=\mathbb{C}[x,y,z]/(f_{E_k})$ for type $E_k$ with $\deg(x)<\deg(y)<\deg(z)$, we see that any graded subalgebra of $R_{E_k}$ generated by elements whose degrees are less than or equal to $d r_{p}(=s r_{p})$ would be contained in $\mathbb{C}[x]$ for any ramification point $p\in \mathrm{Proj}\,R_{E_k}$. This shows that the integral closure of $R/I_S$ cannot be $R_{E_k}$, which is a contradiction. If the integral closure $H^0(\mathcal{O}_{\tilde{S}})$ of $R/I_S$ were isomorphic to the coordinate ring $R_{D_k}=\mathbb{C}[x,y,z]/(f_{D_k})$ for type $D_k$ with $\deg(x)<\deg(y)<\deg(z)$, we can similarly show that $R/I_S$ would be contained in $\mathbb{C}[x,y]$. Then the integral closure of $R/I_S$ cannot be $R_{D_k}$ again.

If $\tilde{S}$ is smooth, we can write $H^0(\mathcal{O}_{\tilde{S}})=\mathbb{C}[u,v]$ as Poisson algebras with $\deg (u)\le \deg(v)$. Since $\mathbb{P}(S)$ contains a point of ramification index $r_{p_5}=d_5/d$ and the ramification indices of $\tilde{S}$ are $\deg (u)/d$ and $\deg(v)/d$, we have $\deg (u)=d_5$ or $\deg(v)=d_5$. If $d_5=\deg (u)< \deg(v)$, then all $\bar{x}_i$'s would be in $\mathbb{C}[u]$ and we have a contradiction. Thus, we have $\deg(v)=d_5$. Since $R/I_S\not\subset \mathbb{C}[u]$ and $p_5\in\mathbb{P}(S)$, we may assume that the homogeneous generator $\bar{x}_5$ is of the form $v+a u^b\,(a\in\mathbb{C},b\ge1)$ by rescaling $x_5$. We show that $R/I_S$ would coincide with $\mathbb{C}[u,v]$, which contradicts with the fact that $S$ is singular. Since the normalization $\tilde{S}\to S$ is a bijection, $R/I_S$ must contain $u^l$ for some $l\ge1$ (see the proof of Lemma \ref{non-normal}). By applying $\{-,\bar{x}_5\}$ to $u^l$ successively, we would obtain $u\in R/I_S$ and we are done.

If $\tilde{S}$ is of type $A_k$, we can write $H^0(\mathcal{O}_{\tilde{S}})=\mathbb{C}[x,y,z]\subset\mathbb{C}[u,v]$ as a Poisson algebra where
$$x=u^{n+1},y=v^{n+1},z=uv$$
and $\deg (u)\le \deg(v)$. Then $\deg (x)$ or $\deg(y)$ is equal to $d_5$ similarly as above. If $d_5=\deg (x)< \deg(y)$, then all $\bar{x}_i$'s would be in $\mathbb{C}[x,z]$, which is a contradiction since in this case the integral closure of $R/I_S$ cannot be $\mathbb{C}[x,y,z]$. Thus, we have $\deg(y)=d_5$. Since the ideal $J\subset R/I_S$ generated by the brackets $\{\bar{x}_i,\bar{x}_j\},1\le i,j\le 5$ becomes the maximal ideal after taking its radical, it must contain $y^l=v^{(n+1)l}$ for some $l\ge1$. For this, some $\bar{x}_{i_0}$ must contain $y^{l'}z=uv^{(n+1)l'+1}$ for some $l'\ge0$ (see the proof of Lemma \ref{non-normal}). However, the condition $d_{i_0}\le d_5=\deg(y)$ implies $l'=0$, and we have $\deg(\bar{x}_{i_0})=s$. Then $S$ is normal by Lemma \ref{non-normal}.
\qed

\vspace{3mm}

From now on, we assume that $p_5\not\in\mathbb{P}(S')$. Since $S'$ is the singular locus of $X$, this is equivalent to the condition that $\partial f/\partial x_i|_{x_1=\cdots=x_4=0}$ is nonzero for some $i$, or in other words $\partial f/\partial x_i$ contains $x_5^k$ for some $i$ and $k\ge1$. Combining with the Pfaffian condition (\ref{4-dim}), the Poisson matrix must satisfy the following condition $(*)_5$:
\[
  \left(
  \begin{tabular}{l}\label{unram5}
there are mutually different $i_1,i_2,i_3,i_4\in\{1,2,\dots,5\}$ such that both $\Theta_{i_1,i_2}$ and \\
$\Theta_{i_3,i_4}$ contain $x_5$.
 \end{tabular}
  \right.
\]
Note that $\Theta_{i,j}$ cannot contain $x_5^2$ since $\deg \Theta_{i,j}=d_i+d_j-s<2d_5$.

The similar statement to Lemma \ref{p5} holds also for $p_4$.

\begin{lem}\label{p4}
Assume $p_5\not\in\mathbb{P}(S')$. If $p_4$ is in $\mathbb{P}(S')$, then $S'$ contains an irreducible component $S$ such that $S$ is a normal symplectic surface of type $A_k$ for some $k\ge1$.
\end{lem}

{\em Proof.} We apply basically the same argument as in the proof of Lemma \ref{p5}. We take the irreducible component $S$ of $S'$ so that $p_4\in \mathbb{P}(S)$ and we define $I_S\subset R$ as before. We first assume that $p_4$ is not a ramification point of $\mathbb{P}(S)$.

If $R/I_S$ is generated by 3 elements, then $S$ is normal but singular since $X$ is indecomposable. We show that $S$ is of type $A_k$. We can choose 3 elements from $\bar{x}_1,\dots,\bar{x}_5$ as generators of $R/I_S$ as before, and one of 3 generators must be $\bar{x}_4$. Note that $p_4\in \mathbb{P}(S)$ implies that vanishing of the other two generators does not give vanishing of $\bar{x}_4$ in $R/I_S$. However, from the classification (Proposition \ref{2-dim}), one sees that this condition excludes $E_k$-type by the maximality or the second maximality of $d_4$ among the degrees of the 3 generators. For type $D_k$, only the point $[0\mathbin:1\mathbin:0]$ in Proposition \ref{ram} is possible as $p_4\in \mathbb{P}(S)$, but this is a ramification point and hence contradiction. Thus we have shown that $S$ is of type $A_k$.

If $R/I_S$ is generated by at least 4 elements, then, in order for $p_4\in \mathbb{P}(S)$ not to be a ramification point, there are two possibilities: (1) four of $d_1,\dots,d_5$ are equal to $d_4$, and (2) $d_2=d_3=d_4$ and $d_5=l d_4$ for some $l\ge2$. By Condition $(*)_5$, we see that (2) does not happen for a degree reason and that (1) falls into two cases:\\
(1-1) $d_4=s$, and\\
(1-2) $d_1=d_2=d_3=d_4\ne s$ and $d_5=2d_4-s$.\\
We can exclude the case (1-1) as in the proof of the previous lemma. In the case (1-2), we have $\Theta_{i,j}\in\mathbb{C}[x_5]$ for $1\le i,j\le 4$ and thus $x_5|f$ by (\ref{4-dim}), which is contrary to the irreducibility of $f$.

Next we assume that $p_4$ is a ramification point of $\mathbb{P}(S)$. We first show that the possible types of the normalization $\tilde{S}$ are $A_k$ and $D_k$. Note that $R/I_S$ is generated by elements whose degrees are less than or equal to $d_4=s r_{p_4}$ and one additional element $\bar{x}_5$ of higher degree. If the integral closure $H^0(\mathcal{O}_{\tilde{S}})$ of $R/I_S$ were isomorphic to the coordinate ring $R_{E_k}=\mathbb{C}[x,y,z]/(f_{E_k})$ for type $E_k$, we see that $R/I_S$ would be contained in $\mathbb{C}[x,\bar{x}_5]$. Then the integral closure of $R/I_S$ cannot be $R_{E_k}$.  We can also exclude the smooth case. Indeed, using the same argument as in the proof of the previous lemma, we may assume that $\bar{x}_4=u$ or $v+a u^b\,(a\in\mathbb{C},b\ge1)$ by rescaling $x_4$, and $R/I_S$ would necessarily be normal as before, which is contrary to the indecomposability of $X$.

We consider the case when $\tilde{S}$ is of type $A_k$ and show that $S=\tilde{S}$. In this case we identify $H^0(\mathcal{O}_{\tilde{S}})$ with $\mathbb{C}[x,y,z]\subset\mathbb{C}[u,v]$ as before and we have that $\deg(x)$ or $\deg(y)$ is equal to $d_4$. If $d_4=\deg(x)<\deg(y)$, then we would have $\bar{x}_1,\dots,\bar{x}_4\in \mathbb{C}[x,z]$. However, $\mathbb{P}(\tilde{S})\cap \{x=z=0\}\ne\emptyset$ implies that $p_5\in \mathbb{P}(S)$, which is a contradiction. If $\deg(x)\le\deg(y)=d_4$, then we may assume that $p_4$ is the ramification point of $\mathbb{P}(\tilde{S})$ defined by $x=z=0$ by switching $x$ and $y$ if necessary (when $\deg(x)=\deg(y)$). This shows that $\bar{x}_1,\bar{x}_2,\bar{x}_3\in \mathbb{C}[x,z]$ and that $\bar{x}_4$ contains $y$. Then some $\bar{x}_{i_0}$ must contain a term $y^{l}z$ for some $l\ge1$ similarly to the proof of the previous lemma. If $l=0$, then we have $\deg(\bar{x}_{i_0})=s$, and hence $S$ is normal by Lemma \ref{non-normal}. If $l\ne0$, then $i_0=5$ and $l=1$ by Condition $(*)_5$ since
$$d_{i_0}\le d_5\le d_3+d_4-s< 2d_4=\deg(y^2).$$
Condition $(*)_5$ also implies that $\bar{\Theta}_{i,j}$ and $\bar{\Theta}_{i',j'}$ for some $\{i,j,i',j'\}=\{1,\dots,4\}$ must contain $yz$, but this is impossible since $\bar{x}_i\in \mathbb{C}[x,z]$ for $1\le i\le3$.

Finally we assume that $\tilde{S}$ is of type $D_k$. Then $H^0(\mathcal{O}_{\tilde{S}})$ is isomorphic to $R_{D_k}=\mathbb{C}[x,y,z]/(f_{D_k})$ with $f_{D_k}=x^{k-1}+x^2y+z^2$ as a graded Poisson algebra where
$$\deg(x)=2s,\deg(y)=(k-2)s\text{ and }\deg(z)=(k-1)s.$$
The ramification point $p_4$ corresponds to $[0\mathbin:1\mathbin:0]$. Indeed, if $r_{p_4}=2$ when $k\ge5$, then $R/I_S$ would be contained in $\mathbb{C}[x,\bar{x}_5]$ similarly for type $E$ above, and the integral closure of $R/I_S$ cannot be $R_{D_k}$. Note that, when $k=4$, the 3 ramification points are symmetric. Thus, $\bar{x}_4\in R/I_S\subset H^0(\mathcal{O}_{\tilde{S}})$ contains $y$ under this identification. Then we also see that $\bar{x}_1,\bar{x}_2$ and $\bar{x}_3$ are 0 or powers of $x$. In particular we have $\{\bar{x}_i,\bar{x}_j\}=0$ for $1\le i,j\le3$. Then, by Condition $(*)_5$, there exist $i\in\{1,2,3\}$ and $j$ such that $\{x_i,x_j\}$ contains $x_5$ and $\{\bar{x}_i,\bar{x}_j\}=0$. Note that, if $\{x_i,x_5\}$ contains $x_5$, then $\bar{x}_i=0$ since $d_i=s$ and $H^0(\mathcal{O}_{\tilde{S}})$ does not contain any nonzero element of degree $s$. Thus, we can write $\bar{x}_5=r(\bar{x}_1,\dots,\bar{x}_4)$ where $r(\bar{x}_1,\dots,\bar{x}_4)$ is a polynomial of $\bar{x}_1,\dots,\bar{x}_4$. This is a contradiction since, in such a case, $R/I_S$ would not contain
$$\{\bar{x}_i,\bar{x}_4\}=\{x^{k_i},y\}=k_ix^{k_i-1}z\;(k_i\ge1)$$
where we take the minimal $i\le3$ with nonzero $\bar{x}_i$, and thus $R/I_S$ would not be closed under Poisson bracket. Note that at least one of $\bar{x}_1,\,\bar{x}_2$, and $\bar{x}_3$ is nonzero since otherwise $R/I_S$ would be generated by two elements.
\qed

\vspace{3mm}

From now on, we also assume that $p_4\not\in\mathbb{P}(S')$. Similarly to the condition $p_5\not\in\mathbb{P}(S')$, the partial derivative $\partial f/\partial x_i$ contains $x_4^k$ for some $i$ and $k\ge2$. By considering possible degrees $d_i$ using Condition $(*)_5$, one can show that the Poisson matrix must satisfy the following condition $(*)_4$:
\[
  \left(
\begin{tabular}{l}\label{unram4}
There are mutually different $i_1,i_2,i_3,i_4\in\{1,2,\dots,5\}$ such that $\Theta_{i_1,i_2}$ contains $x_4$\\
and $\Theta_{i_3,i_4}$ contains $x_4^{l_0}$ for $l_0=1$ or $2$.
 \end{tabular}
  \right.
\]
Indeed, we have $\deg(\Theta_{i,j})=d_i+d_j-s<2d_4$ if $i,j\ne5$, and, for any $i\le 4$, we also have
$$\deg(\Theta_{i,5})=d_i+d_5-s=d_i+(d_j+d_k-s)<3d_4$$
for some $j,k\le4$ by Condition $(*)_5$.

\begin{lem}\label{*45}
Assume that Condition $(*)_5$ holds. If $d_i=s$ for some $i$, then any irreducible component $S$ of $S'$ is normal of type $A_k$ for some $k\ge1$.
\end{lem}

{\em Proof.} If $d_i=s$ for some $i$, then there exists $i_0\ne5$ such that $d_{i_0}=s$ and $\Theta_{i_0,j}$ contains $x_5$ for some $j$ by Condition $(*)_5$. If $\bar{x}_{i_0}\ne0\in R/I_S$, then $S$ is normal of type $A_k$ for some $k\ge1$ by Lemma \ref{non-normal}. So we assume $\bar{x}_{i_0}=0$. Then we have $\bar{\Theta}_{i_0,j}=0$, and this implies that $\bar{x}_5$ is a polynomial of $\bar{x}_1,\dots,\bar{x}_4\in R/I_S$. Therefore, $R/I_S$ is generated by 3 elements and hence integrally closed. We assume that $S$ were not of type $A_k$ to deduce contradiction.

We have $i_0=1$ since $R/I_S$ for $D_k$ or $E_k$ is generated by homogeneous elements of degree greater than $s$. By replacing generators, we may assume that $I_S=(x_1,x_5,\Delta)$ where $\Delta\in\mathbb{C}[x_2,x_3,x_4]$ is the equation for $D_k$ or $E_k$, and we can write the defining equation of $X$ as
$$f=h(x_2,x_3,x_4)\Delta^2+P(x_1,\dots,x_5)$$
(see the argument before (\ref{f}) in Subsection \ref{4.2}). Thus, we obtain $\deg(f)\ge 2\deg(\Delta)$. We may assume that Condition $(*)_5$ still holds after the coordinate change by Lemma \ref{p5}. In particular we have $d_3+d_4-s\ge d_5$. Since $\deg(f)=d_1+\cdots+d_5-2s$ by (\ref{Pf}) and $\deg(\Delta)=d_2+d_3+d_4-s$, we have
$$\begin{aligned}
0&\le \deg(f)-2\deg(\Delta)=d_1+d_5-2s-(d_2+d_3+d_4-2s)\\
&=d_5+s-(d_2+d_3+d_4)\le d_5+s-(d_2+d_5+s)=-d_2.
\end{aligned}$$
This is a contradiction since $d_2>0$.
\qed

\vspace{3mm}

If Condition $(*)_5$ and $(*)_4$ for $l_0=1$ are satisfied, then one can easily check that there exists $i$ with $d_i=s$. Therefore, in his case, $X$ contains a normal surface of type $A_k$ by Lemma \ref{*45}. So we assume that Condition $(*)_5$ and $(*)_4$ for $l_0=2$ are satisfied from now on. We may also assume that there are no $i$ with $d_i=s$ by Lemma \ref{*45}. Then we see that we have the following 3 possibilities:

\vspace{3mm}

\noindent
\underline{Case 1.} Both $\Theta_{1,4}$ and $\Theta_{2,3}$ contain $x_5$. $\Theta_{1,2}$ contains $x_4$ and $\Theta_{3,5}$ contains $x_4^2$.
$$(d_1,\dots,d_5)=(a,2a-2s,2a-s,3a-3s,4a-4s)$$
for some $a\ge2s$. In this case we may assume that $\Theta_{1,2}=x_4$ and $\Theta_{1,4}=x_5$ by replacing the coordinates $x_4$ and $x_5$.

\vspace{3mm}

\noindent
\underline{Case 2.} Both $\Theta_{1,4}$ and $\Theta_{2,3}$ contain $x_5$. $\Theta_{1,2}$ contains $x_4$ and $\Theta_{4,5}$ contains $x_4^2$.
$$(d_1,\dots,d_5)=(2s,a,3s,a+s,a+2s)$$
for some $2s\le a\le3s$. In this case we may assume that $\Theta_{1,2}=x_4$ and $\Theta_{1,4}=x_5$ by replacing the coordinates $x_4$ and $x_5$.

\vspace{3mm}

\noindent
\underline{Case 3.} Both $\Theta_{1,4}$ and $\Theta_{2,3}$ contain $x_5$. $\Theta_{1,3}$ contains $x_4$ and $\Theta_{4,5}$ contains $x_4^2$.
$$(d_1,\dots,d_5)=(2s,3s,a,a+s,a+2s)$$
for some $a\ge3s$. In this case we may assume that $\Theta_{1,3}=x_4$ and $\Theta_{1,4}=x_5$ by replacing the coordinates $x_4$ and $x_5$.

\vspace{3mm}

We can show that these cases do not occur except $a=3s$ in Case 1.

\begin{lem}\label{p3}
Assume that $p_i\not\in\mathbb{P}(S')$ for $i=4,5$ and that $d_i\ne s$ for any $i$. Then we have $(d_1,\dots,d_5)=(3s,4s,5s,6s,8s)$.
\end{lem}

{\em Proof.} We treat Case 1 and Case 2,3 separately. For Case 1, we show that $p_3\not\in\mathbb{P}(S)$, and then show that the only case $a=3s$ is possible. We assume $p_3\in\mathbb{P}(S)$ for an irreducible component $S$ of $S'$ to deduce contradiction. Then $p_3$ is a ramification point of $\mathbb{P}(S)$. Indeed, otherwise at least three of $d_1,\dots,d_5$ would be multiples of $d_3$, which is impossible. 

If the normalization $\tilde{S}$ of $S$ were of type $E_k$, then the possible ramification index $r_{p_3}=d_3/s$ of $p_3$ would be $3,4$ or $5$ since $d_3\ge 3s$ (see Proposition \ref{ram}). But in any case the coordinate ring for $E_k$ does not contain elements of degree smaller than $d_3=sr_{p_3}$, which implies that $\bar{x}_1=\bar{x}_2=0$ and hence $\bar{x}_4=\bar{\Theta}_{1,2}=0$. This is a contradiction.

If $\tilde{S}$ were of type $D_k$, then the possible ramification index $r_{p_3}=d_3/s$ of $p_3$ would be $k-2$ with $k\ge5$ since $d_3\ge3s$. Then we would have $\bar{x}_1,\bar{x}_2\in\mathbb{C}[x]$ when we identify $R/I_S$ with the coordinate ring $R_{D_k}=\mathbb{C}[x,y,z]/(f_{D_k})$ of type $D_k$ with $\deg(x)=2s$. This would give $\bar{x}_4=\bar{\Theta}_{1,2}=0$ and hence $\bar{x}_5=\bar{\Theta}_{1,4}=0$, which is a contradiction since we would have $R/I_S\subset \mathbb{C}[x,y]$ and thus the integral closure of $R/I_S$ cannot be $R_{D_k}$.

If $\tilde{S}$ were of type $A_k$ or smooth, then the same argument as in the proof of Lemma \ref{p5} (or \ref{p4}) shows that we may assume that $\bar{x}_3$ is of the form $x+p_1(z)$ or $y+p_2(x,z)$ for some polynomial $p_i$ when we identify $R/I_S$ with the coordinate ring $\mathbb{C}[x,y,z]\subset \mathbb{C}[u,v]$ (possibly with $k=0$) as before. If $k=0$, then  the same argument as in  the proof of Lemma \ref{p5} shows that $S$ is smooth, which is a contradiction. If $k\ne0$, then some $\bar{x}_{i_0}$ must contain $y^l z$ for some $l\ge0$ by the same argument in Lemma \ref{p5}. This $l$ is not equal to $0$ since $d_i\ne s=\deg(z)$ for any $i$. Since $d_2<d_3\le\deg(y)$, we would have $\bar{x}_1,\bar{x}_2\in\mathbb{C}[x,z]$ and thus $\bar{x}_4=\bar{\Theta}_{1,2}$ and $\bar{x}_5=\bar{\Theta}_{1,4}$ are also in $\mathbb{C}[x,z]$. This is contrary to the existence of $\bar{x}_{i_0}$ which contains $y^l z$.

Thus, we have shown that $p_3\not\in\mathbb{P}(S)$. Then the Poisson structure on $\mathbb{C}[x_1,\dots,x_5]$ must satisfy the condition $(*)_3$ which is similarly defined to $(*)_5$ and $(*)_4$. However, one can check that this condition cannot be satisfied unless $a=3s$ by considering the possible degrees of $\Theta_{i,j}$'s. Note that in the case $a=3s$, $\Theta_{1,5}$ and $\Theta_{3,4}$ can contain $x_3^2$. This completes the Case 1.

Next we show that Case 2 and 3 are impossible. In order to treat these two cases at once, we switch $x_2$ and $x_3$ for Case 2 (We will not use the inequality $d_2\le d_3$ in the proof). Then the defining equation $f$ of $X$ contains a monomial $x_3 x_4^3$ because of the condition $(*)_4$. Since we have $\Theta_{1,3}=x_4$, the bracket $\{x_1,x_3 x_4^3\}$ must contain $x_4^4$. We may assume that $\Theta_{1,2}$ does not contain $x_4$ (when $d_2=d_3$) by replacing the coordinate $x_2$ with $x_2+cx_3$ for some $c\in\mathbb{C}$. 
Then one can easily check, by using the Leibniz rule, that $x_3 x_4^3$ is a unique monomial $g\in\mathbb{C}[x_1,\dots,x_5]$ such that $\{x_1,g\}$ contains $x_4^4$. However, this implies that $\{x_1,f\}$ contains $x_4^4$, which is contrary to the condition $\{x_1,f\}=0$ (\ref{Pf}).
\qed

\vspace{5mm}

Therefore, we are reduced to the following exceptional case:
$$(d_1,\dots,d_5)=(3,4,5,6,8)$$
where we set $s=1$ so that $\mathrm{g.c.d.}(d_1\dots,d_5)=1$. Recall that we may assume that $\Theta_{1,2}=x_4,\;\Theta_{1,4}=x_5,$ and $\Theta_{2,3}=x_5+\cdots$ by suitably replacing the coordinates $x_4$ and $x_5$ and by rescaling $x_3$. We show that this case cannot happen by a direct calculation (using a computer). The computation is carried out as follows. 
We consider the possible values for all $\Theta_{i,j}$ so that $\Theta_{i,j}$'s are homogeneous:
$$\begin{aligned}
&\Theta_{1,2}=x_4,\;\Theta_{1,3}=a_1 x_1 x_2,\;\Theta_{1,4}=x_5,\;\Theta_{1,5}=a_2 x_1^2 x_4+a_3 x_2 x_4+a_4 x_3^2,\\
&\Theta_{2,3}=x_5+a_5 x_1 x_3+a_6 x_2^2,\Theta_{2,4}=a_7 x_1 x_4+a_8 x_2 x_3+a_9 x_1^3,\\
&\Theta_{2,5}=a_{10}x_1 x_5+a_{11} x_3 x_4+a_{12} x_1^2 x_3+a_{13} x_1 x_2^2,\;\Theta_{3,4}=a_{14} x_1^2 x_2+a_{15}x_3^2+a_{16}x_2 x_4,\\
&\Theta_{3,5}=a_{17}x_2 x_5+a_{18} x_4^2+a_{19} x_1^2 x_4+a_{20} x_1 x_2 x_3+a_{21} x_2^3+a_{22} x_1^4,\\
&\Theta_{4,5}=a_{23} x_3 x_5+a_{24} x_1 x_2 x_4+a_{25} x_1 x_3^2+a_{26} x_2^2 x_3+a_{27} x_1^3 x_2
\end{aligned}$$ 
using parameters $a_1,a_2,\cdots,a_{27}$. Then we compute the Jacobi relations
$$J_{i,j,k}:=\{x_i,\{x_j,x_k\}\}+\{x_j,\{x_k,x_i\}\}+\{x_k,\{x_i,x_j\}\}$$
for all $1\le i<j<k\le5$. The following formula is useful in calculation
$$\{x_i,g\}=\sum_{j=1}^5 \partial g/\partial x_j \{x_i,x_j\},\;g\in\mathbb{C}[x_1,\dots,x_5].$$
Each $J_{i,j,k}$ is a polynomial in $x_1,\dots,x_5$ with coefficients in $\mathbb{C}[a_1,a_2,\dots]$. All these coefficients in $J_{i,j,k}$'s must vanish in order for the Poisson bracket to satisfy Jacobi identities. We compute (the radical of) the ideal generated by the coefficients in $J_{i,j,k}$'s. As in the proof of Lemma \ref{p3}, we may assume that the coefficients $a_4$ and $a_{15}$ of $x_3^2$ must be nonzero. We see that the solution space $V\subset \mathbb{C}^{27}$ has a unique irreducible component $V'$ which is not contained in $\{a_4 a_{15}=0\}$. The component $V'$ is an affine line defined by
$$\begin{aligned}
a_1&=a_2=a_3=a_5-4a_4=a_6=a_7+2a_4=a_8-4a_4=a_9+20a_4=a_{10}+2a_4=a_{11}-4a_4\\
&=a_{12}=a_{13}=a_{14}=a_{15}+a_4=a_{16}=a_{17}=a_{18}=a_{19}=a_{20}=a_{21}=a_{22}=a_{23}+2a_4\\
&=a_{24}=a_{25}+5a_4=a_{26}=a_{27}=0.
\end{aligned}$$ 
In particular the coefficient $a_{18}$ of $x_4^2$ vanishes on $V'$. This is contrary to the assumption on Case 1. Therefore, we have shown that $X$ contains a type-$A$ normal symplectic surface in all possible cases, and Proposition \ref{A} is proved.


\subsection{$X$ is isomorphic to $X_n$}\label{4.2}

In the previous subsection we have shown that $S'=\mathrm{Sing}(X)$ contains an irreducible component $S$ which is a conical symplectic surface of type $A_k$. In this subsection we complete the proof of Theorem \ref{main} using this fact. More concretely, we show that the Poisson matrix $\Theta$ associated to $X$ can be transformed into $\Theta_n$ associated to $X_n$ (see Section \ref{2}) by using suitable coordinates of $\mathbb{C}^5$.

The key point is that the graded polynomial ring $\mathbb{C}[x_1,\dots,x_5]$ has a distinguished generator, which we will denote by $x_{\alpha_1}$, of the same degree as the weight of the symplectic form due to the existence of the normal surface of type $A_k$. The linear operator obtained by taking the Poisson bracket with this distinguished element gives a new grading $w(-)$. The important feature of this grading is that Poisson brackets of homogeneous elements are again homogeneous with respect to $w(-)$ (Lemma \ref{wt}). 

\vspace{4mm}

{\em Proof of Theorem \ref{main}}\\
We first make the ideal of the surface $S$ into a normal form. Since the dimension of the tangent space of $S$ at the origin is $3$, the defining ideal $I\subset \mathbb{C}[x_1,\dots,x_5]$ of $S$ in $\mathbb{C}^5$ must contain two elements $f_1$ and $f_2$ whose linear terms are linearly independent. By applying a linear coordinate change, we may assume that the linear terms of $f_1$ and $f_2$ are $x_{\alpha_4}$ and $x_{\alpha_5}$ respectively where ${\alpha_4},{\alpha_5}\in\{1,\dots,5\}$ with ${\alpha_4}< {\alpha_5}$. Since $x_i$'s have positive degrees, $f_1-x_{\alpha_4}$ consists of variables different from $x_{\alpha_4}$ and $x_{\alpha_5}$. Note that we have $d_1\le \cdots\le d_5$ by assumption. Thus,
$$(\{x_1,\dots,x_5\}\setminus \{x_{\alpha_4}\})\cup \{f_1\}$$
is a new generating system of $\mathbb{C}[x_1,x_2,\dots,x_5]$. Similarly, $f_2-x_{\alpha_5}$ consists of variables different from $x_{\alpha_5}$, and
$$(\{x_1,\dots,x_5\}\setminus \{x_{\alpha_4},x_{\alpha_5}\})\cup \{f_1,f_2\}$$
is a new generating system of $\mathbb{C}[x_1,x_2,\dots,x_5]$. Therefore, we may assume that $S$ is defined by three functions $x_{\alpha_4},x_{\alpha_5}$ and $\Delta\in\mathbb{C}[x_{\alpha_1},x_{\alpha_2},x_{\alpha_3}]$ by replacing $x_{\alpha_4}$ and $x_{\alpha_5}$ by $f_1$ and $f_2$ respectively where $\{{\alpha_1},{\alpha_2},{\alpha_3},{\alpha_4},{\alpha_5}\}=\{1,2,\dots,5\}$ and $\Delta$ is a defining equation of $S$ in $\mathbb{C}^3$. By  Proposition \ref{2-dim}, we may assume that
$$\Delta=\frac{1}{k+1}x_{\alpha_1}^{k+1}+x_{\alpha_2} x_{\alpha_3}$$
by replacing generators of $\mathbb{C}[x_{\alpha_1},x_{\alpha_2},x_{\alpha_3}]$. Then we have $d_{\alpha_1}=s$, and we may assume that
$\alpha_2< \alpha_3$ and $\alpha_4<\alpha_5$.

Since the defining equation $f$ of $X$ is in $I=(x_{\alpha_4},x_{\alpha_5},\Delta)$, it is of the form
$$f=h_0(x_{\alpha_1},x_{\alpha_2},x_{\alpha_3})\Delta+p(x_1,x_2,x_3,x_4,x_5)$$
where $h_0(x_{\alpha_1},x_{\alpha_2},x_{\alpha_3})$ is a nonzero polynomial and $p(x_1,x_2,x_3,x_4,x_5)$ is in the ideal $(x_{\alpha_4},x_{\alpha_5})$. Since $\partial f/\partial x_{\alpha_i}$ vanishes on $S=\{x_{\alpha_4}=x_{\alpha_5}=\Delta=0\}$ for $i=1,2,3$, the polynomial $h_0(x_{\alpha_1},x_{\alpha_2},x_{\alpha_3})$ is divisible by $\Delta$: 
\begin{equation}\label{f}
f=h(x_{\alpha_1},x_{\alpha_2},x_{\alpha_3})\Delta^2+p(x_1,x_2,x_3,x_4,x_5).
\end{equation}

We denote by $\bar{g}$ the image in $R/I$ of an element $g\in \mathbb{C}[x_1,x_2,\dots,x_5]$. Then we have $\bar{x}_i=0$ for $i={\alpha_4},{\alpha_5}$, and also have
$$\bar{\Theta}_{{\alpha_1},{\alpha_2}}=\bar{x}_2,\;\bar{\Theta}_{{\alpha_1},{\alpha_3}}=-\bar{x}_3,\text{ and  }\;\bar{\Theta}_{{\alpha_2},{\alpha_3}}=\bar{x}_1^k.$$
Note that $I$ is a Poisson ideal since $S$ is an irreducible component of the singular locus of $X$. In particular, we have $\bar{\Theta}_{i,\alpha_4}=\bar{\Theta}_{i,\alpha_5}=0$ for any $i$.

We can write
$$\Theta_{{\alpha_1},{\alpha_4}}=a_1 x_{\alpha_4}+p_{14},\;\Theta_{{\alpha_1},{\alpha_5}}=a_2 x_{\alpha_5}+p_{15}$$
where $a_1,a_2\in\mathbb{C},\, p_{14}\in (x_{\alpha_5},\Delta)$ and $p_{15}\in (x_{\alpha_4},\Delta)$. Then we define complex numbers $w_i\,(i=1,\dots,5)$ as $$w_{\alpha_1}=0,\,w_{\alpha_2}=1,\,w_{\alpha_3}=-1,\,w_{\alpha_4}=a_1,\text{ and }\;w_{\alpha_5}=a_2.$$
For a nonzero monomial $g=ax_1^{k_1}\cdots x_5^{k_5}$, we define $w(g)$ as $\sum_{i=1}^5 k_i w_i$. If a polynomial $g=\sum_j g_j\in\mathbb{C}[x_1,x_2,\dots,x_5]$ consists of monomials $g_j$ such that $w(g_j)$ are the same for all $j$, then we say that $g$ is $w$-{\em homogeneous} of degree $w(g):=w(g_j)$. This defines a $\mathbb{C}$-grading on $\mathbb{C}[x_1,x_2,\dots,x_5]$ which is different from the one defined from the conical $\mathbb{C}^*$-action on $X$.

\begin{lem}\label{wt}
There are homogeneous coordinates $x'_1,x'_2,\dots,x'_5$ of $\mathbb{C}[x_1,x_2,\dots,x_5]$\\
(i.e., $\mathbb{C}[x'_1,x'_2,\dots,x'_5]=\mathbb{C}[x_1,x_2,\dots,x_5]$) and a permutation $\beta_1,\beta_2,\dots,{\beta_5}$ of $1,2,\dots,5$ such that\\
$\bullet\;\{x'_{\beta_1},x'_i\}=w'_i x'_i+p_i(x'_1,\dots, x'_{i-1})$ for some $w'_i\in\mathbb{C}$ and $p_i\in\mathbb{C}[x'_1,\dots,x'_{i-1}]$ $(\forall i)$,\\
$\bullet$ $\{x'_{\beta_1},x'_i\}$ is $w'$-homogeneous of degree $w'_i\,(\forall i)$,\\
$\bullet\;\deg(x'_i)=d_i\,(\forall i)$, and\\
$\bullet\;I=(x'_{\beta_4},x'_{\beta_5},\Delta')$ with $\Delta'=\frac{1}{k+1} {x'_{\beta_1}}^{k+1}+x'_{\beta_2} x'_{\beta_3}$.\\
Here, the grading $w'(-)$ is defined by setting $w'(x'_i)=w'_i$ (similarly to $w(-)$). Moreover, $\{x'_i,x'_j\}$ is $w'$-homogeneous of degree $w'_i+w'_j$ for all nonzero $\{x'_i,x'_j\}$.
\end{lem}

{\em Proof.} For each $i$, we can write
$$\Theta_{{\alpha_1},i}=w_i x_i +q_i(x_1,\dots, x_{i-1})+b_{i,1} x_{\alpha_4}+b_{i,2} x_{\alpha_5}$$
with $q_i\in I$ and $b_{i,j}\in\mathbb{C}$ satisfying $b_{{\alpha_4},1}=0$ and $b_{{\alpha_5},1}=b_{{\alpha_5},2}=0$. Note that $I$ does not contain $x_{\alpha_1},\,x_{\alpha_2}$ and $x_{\alpha_3}$, and thus $\Theta_{{\alpha_1},i}-w_i x_i$, which is in $I$, does not contain linear terms of $x_{\alpha_1},\,x_{\alpha_2}$ and $x_{\alpha_3}$ as well. Note also that $\alpha_4<\alpha_5$ by assumption and in particular the linear term $b x_{\alpha_4}$ in $\Theta_{{\alpha_1},\alpha_5}$ can be taken in $q_{\alpha_5}$.

We can assume that $b_{i,1}=b_{i,2}=0$ for all $i$ by a suitable linear change of coordinates. This is shown as follows. For $i={\alpha_1}$, this always holds since $\Theta_{{\alpha_1},{\alpha_1}}=0$. For $i={\alpha_4}$ and ${\alpha_5}$, let $$\rho:\mathbb{C}x_4\oplus\mathbb{C}x_5\to\mathbb{C}x_4\oplus\mathbb{C}x_5$$
be the $\mathbb{C}$-linear map obtained by restricting $$\{x_{\alpha_1},-\}:\mathbb{C}[x_1,x_2,\dots,x_5]\to\mathbb{C}[x_1,x_2,\dots,x_5]$$
to $\mathbb{C}x_4\oplus\mathbb{C}x_5$ and extracting just linear terms from the image. Then we can write
$$\rho(x_{\alpha_4})=w_{\alpha_4}x_{\alpha_4}+b_{{\alpha_5},1}x_{\alpha_5},\,\rho(x_{\alpha_5})=w_{\alpha_5}x_{\alpha_5}+bx_{\alpha_4}.$$
Then we can make the matrix representation of $\rho$ into an upper triangular one (i.e. $b_{{\alpha_5},1}=0$) by choosing a suitable basis of $\mathbb{C}x_4\oplus\mathbb{C}x_5$ by linear algebra. Finally, for $i=\alpha_2$ and $\alpha_3$, we can assume that $b_{i,1}=b_{i,2}=0$ just by reordering $x_i$'s when $b_{i,j}$ is nonzero (and hence $d_i=d_{\alpha_4}$ or $d_{\alpha_5}$) so that the linear terms $b_{i,1} x_{\alpha_4}$ and $b_{i,2} x_{\alpha_5}$ can be taken in $q_i$. Thus, we may assume that $b_{i,1}=b_{i,2}=0$ for all $i$ and the first property in the statement of the lemma is satisfied. Note that the above coordinate changes can be done preserving the degrees of $x_i$ and that the condition $I=(x_{\beta_4},x_{\beta_5},\Delta'')$ holds by choosing $\beta_i$'s suitably (according to the reordering of $x_i$'s) where $\Delta''=\frac{1}{k+1}x_{\beta_1}^{k+1}+x_{\beta_2} x_{\beta_3}$. From now on we use the coordinate system $x_{\beta_1},\dots,x_{\beta_5}$ obtained in this way.

We define a total order $\prec$ on the set of (monic) monomials of $\mathbb{C}[x_1,x_2,\dots,x_5]$ as follows:
$$x_1^{k_1}\cdots x_5^{k_5}\prec x_1^{l_1}\cdots x_5^{l_5}\underset{\text{def}}{\iff} \text{there exists }i_0\text{ such that }k_{i_0}<l_{i_0}\text{ and }\;k_i=l_i\text{ for }i>i_0.$$
We set $x'_{\beta_1}:=x_{\beta_1}$. For the second property in the statement of the lemma, we use induction on $i$. The case $i=1$ is clear since $q_1=0$. We assume that the claim holds up to $i-1$ for suitable coordinates $x'_1,x'_2,\dots,x'_{i-1},x_i,\dots,x_5$. We show that the claim holds also for $i$ by replacing $x_i$ by another coordinate $x_i'$. If the claim does not hold with respect to $x_i$, we can take a nonzero term $q_i'$ of $q_i$ such that $w'(q_i')\ne w_i=w'(x_i)$. We choose the maximal $q_i'$ with respect to $\prec$. Then replacing $x_i$ by $x_i+\frac{1}{w_i-w'(q_i')}q_i'$ gives new $q_i$ which consists of smaller terms than $q_i'$ with respect to $\prec$. Indeed, we have
$$\begin{aligned}
\left\{x'_{\beta_1}, x_i+\frac{1}{w_i-w'(q_i')}q_i'\right\}&=\{x'_{\beta_1}, x_i\}+\left\{x'_{\beta_1}, \frac{1}{w_i-w'(q_i')}q_i'\right\}\\
&=(w_i x_i +q_i'+r_i)+\left(\frac{1}{w_i-w'(q_i')}(w'(q_i')q_i'+r_i')\right)\\
&=w_i\cdot\left(x_i+\frac{1}{w_i-w'(q_i')}q_i'\right)+r_i+\frac{1}{w_i-w'(q_i')}r_i'
\end{aligned}$$
where $r_i=q_i-q_i'$ and $r_i'=\{x'_{\beta_1},q_i'\}-w'(q_i')q_i'$. Then we can show that the terms of $r_i'$ are smaller than $q_i'$ with respect to $\prec$ by combining the Leibniz rule and the fact that $\{x'_{\beta_1},x'_j\}-w'_j x_j$ is in $\mathbb{C}[x'_1,\dots, x'_{j-1}]$ for $j<i$. Note that any monomial of $\mathbb{C}[x'_1,\dots, x'_{j-1}]$ is smaller than $x'_j$ with respect to $\prec$.

We iterate this process until $\{x'_{\beta_1},x_i\}$ becomes $w'$-homogeneous of degree $w_i$, and we define $x'_i$ to be the resulting coordinate. Since each $d_i$ is positive, possible monomials in $q_i$ are finite and this process will stop in finite time. This process does not preserve the condition $I=(x_{\beta_4},x_{\beta_5},\Delta'')$ in general since $q_i'$ might not be in $I$. In the above process, however, we may use $g\Delta''$ instead of $q_i'$ for a suitable $w'$-homogeneous $g\in\mathbb{C}[x_{\beta_1},x_{\beta_2},x_{\beta_3}]$ with $w'(g\Delta'')=w'(q_i')$ if $q_i'$ is not in $(x_{\beta_4},x_{\beta_5})$. (This happens only when $i=\beta_4$ or $\beta_5$ since $\deg(\Delta'')>d_i$ for $i=\beta_1,\beta_2,\beta_3$.) Indeed, the sum of the terms of $q_i$ which are not in $(x_{\beta_4},x_{\beta_5})$ is of the form $g_1 \Delta''$ where $g_1\in\mathbb{C}[x_{\beta_1},x_{\beta_2},x_{\beta_3}]$. If we choose $g$ as the maximal term of $g_1$ with respect to $\prec$, then $g\Delta''$ contains $q_i'$ as the maximal term. Then we can carry out the procedure using $g \Delta''$ since $\Delta''$ is $w$-homogeneous (with $w(\Delta'')=0$). Now the procedure preserves the condition $I=(x_{\beta_4},x_{\beta_5},\Delta'')$, and thus we have $I=(x'_{\beta_4},x'_{\beta_5},\Delta')$.

For the final claim, we use induction on $i+j$. By the Jacobi identity $J_{\beta_1, i, j}=0$, we have
\begin{equation}\label{i+j}
\{x'_{\beta_1},\{x'_i,x'_j\}\}=(w'_i+w'_j)\{x'_i,x'_j\}+\{p_i,x'_j\}+\{x'_i,p_j\}.
\end{equation}
Note that $\{x'_{\beta_1},g\}$ is $w'$-homogeneous of degree $w'(g)$ for any monomial $g\in\mathbb{C}[x'_1,x'_2,\dots,x'_5]$ since this holds for generators $g=x'_i$ for any $i$ and we can extend it to any polynomial by the Leibniz rule. Note also that $\{p_i,x'_j\}+\{x'_i,p_j\}$ is $w'$-homogeneous of degree $w'_i+w'_j$ by induction hypothesis and the Leibniz rule since $p_i\in\mathbb{C}[x'_1,\dots, x'_{i-1}]$. Let $q$ be a nonzero monomial of $\{x'_i,x'_j\}$ which is maximal with respect to $\prec$. If $q$ is a monomial of $\{p_i,x_j\}+\{x_i,p_j\}$, then $w'(q)=w'_i+w'_j$. If $q$ is not contained in $\{p_i,x'_j\}+\{x'_i,p_j\}$, then $q$ appears with the coefficient $w'_i+w'_j$ on the right hand side of (\ref{i+j}). Since $\{x'_{\beta_1},q\}$ is the sum of $w'(q)q$ and smaller monomials than $q$ with respect to $\prec$, the maximality of $q$ implies that the coefficient of $q$ in the left hand side of (\ref{i+j}) is equal to  $w'(q)$, and thus we obtain $w'(q)=w'_i+w'_j$. Let $q'$ be a nonzero next larger monomial of $\{x'_i,x'_j\}$ with respect to $\prec$. If $\{x'_{\beta_1},q\}$ contains $q'$, then $w'(q')=w'(q)=w'_i+w'_j$ by the second claim of this lemma. If not, the same argument as above shows again that $w'(q')=w'_i+w'_j$. By proceeding the same argument for all the monomials of $\{x'_i,x'_j\}$, we see that $\{x'_i,x'_j\}$ is $w'$-homogeneous of degree $w'_i+w'_j$.
\qed

\vspace{3mm}

By using the coordinates $x'_1,x'_2,\dots,x'_5$ in Lemma \ref{wt}, we may assume that $x_1,\dots, x_5$ are $w$-homogeneous and that $\Theta_{i,j}$ is $w$-homogeneous of degree $w_i+w_j$ for all nonzero $\Theta_{i,j}$. We have $w(f)=\sum_{i=1}^5 w_i=a_1+a_2$ by (\ref{4-dim}). We also have
$$\{x_{\alpha_1},f\}=\sum_{i=1}^5 \Theta_{\alpha_1,i} \frac{\partial f}{\partial x_i}=0$$
by (\ref{Pf}). Note that $\mathrm{pf}(\Theta)$ is in the kernel of $\Theta$. Therefore, we have $a_1+a_2=0$.

\vspace{3mm}

\noindent
Case (I): $a_1=-a_2\ne0$\\
We show that $X$ is isomorphic to $X_n$ as a Poisson scheme (but not as a $\mathbb{C}^*$-variety) in this case.

First we show that $\{x_{\alpha_4},\Theta_{{\alpha_2},{\alpha_3}}\}$ contains $x_{\alpha_1}^{k-1} x_{\alpha_4}$. For this, note that we have
$$\deg(\Theta_{{\alpha_2},{\alpha_3}})=d_{\alpha_2}+d_{\alpha_3}-s<d_{\alpha_2}+d_{\alpha_3}=\deg(\Delta).$$
This implies that  the polynomial $\Theta_{{\alpha_2},{\alpha_3}}-x_{\alpha_1}^k\in I$ is in $(x_{\alpha_4},x_{\alpha_5})$. The fact that  $w(\Theta_{{\alpha_2},{\alpha_3}})=w_{\alpha_2}+w_{\alpha_3}=0$ also implies that $\Theta_{{\alpha_2},{\alpha_3}}-x_{\alpha_1}^k$ is in the ideal
$$I'=(x_{\alpha_2}x_{\alpha_4},x_{\alpha_3}x_{\alpha_4},x_{\alpha_2}x_{\alpha_5},x_{\alpha_3}x_{\alpha_5},x_{\alpha_4}x_{\alpha_5})$$
since $a_1,a_2\ne0$. One can check that $\{x_{\alpha_4},p_0\}$ does not contain $x_{\alpha_1}^{k-1}x_{\alpha_4}$ for any $p_0\in I'$ by using the Leibniz rule and the fact that $\{x_{\alpha_4},x_i\}$ does not contain $x_{\alpha_1}^{k-1}$ for any $i$ since $x_{\alpha_1}^{k-1}\not\in I$. Since $\{x_{\alpha_4},x_{\alpha_1}^k\}$ contains $x_{\alpha_1}^{k-1} x_{\alpha_4}$, the claim follows.

The Jacobi identity $J_{\alpha_2, \alpha_3, {\alpha_4}}=0$ shows that $\{x_{\alpha_2},\Theta_{{\alpha_3},{\alpha_4}}\}$ or $\{x_{\alpha_3},\Theta_{{\alpha_2},{\alpha_4}}\}$ contains $x_{\alpha_1}^{k-1}x_{\alpha_4}$. If $\{x_{\alpha_2},\Theta_{{\alpha_3},{\alpha_4}}\}$ contains $x_{\alpha_1}^{k-1}x_{\alpha_4}$, then $\Theta_{{\alpha_3},{\alpha_4}}$ contains $x_{\alpha_1}^l x_{\alpha_5}$ and
$\Theta_{{\alpha_2},{\alpha_5}}$ contains $x_{\alpha_1}^m x_{\alpha_4}$ for some $l,m\ge0$ with $l+m=k-1$. This is shown as follows. Suppose that $\Theta_{{\alpha_3},{\alpha_4}}$ contains
$$\lambda:=x_{\alpha_1}^{l_1}x_{\alpha_2}^{l_2}x_{\alpha_3}^{l_3}x_{\alpha_4}^{l_4}x_{\alpha_5}^{l_5}$$
and that $\{x_{\alpha_2},\lambda\}$ contains $x_{\alpha_1}^{k-1}x_{\alpha_4}$. Then we have $l_3+l_5\le1$ since otherwise every term of $\{x_{\alpha_2},\lambda\}$ would be divided by $x_{\alpha_3}$ or $x_{\alpha_5}$. Similarly, we have $l_2=0$ and $l_4\le 2$. Since $w(\lambda)=w_{\alpha_3}+w_{\alpha_4}=a_1-1$ and $a_1+a_2=0$, possible $(l_3,l_4,l_5)$ are
$$(0,0,0),(0,0,1),\,(0,1,1),\,(0,2,0),\,(1,1,0),\text{ and }\;(1,2,0).$$
The last two cases do not happen since $\deg(\Theta_{{\alpha_3},{\alpha_4}})<\deg(x_{\alpha_3}x_{\alpha_4})$. In the case $(0,0,0)$, $l_1$ would be greater than $k$ since $\Theta_{{\alpha_3},{\alpha_4}}$ is in $I$, but then $\{x_{\alpha_2},x_{\alpha_1}^{l_1}\}$ cannot contain $x_{\alpha_1}^{k-1}x_{\alpha_4}$, a contradiction. The case $(0,1,1)$ and $(0,2,0)$ also do not happen. Indeed, in these cases $\Theta_{{\alpha_2},{\alpha_5}}$ and $\Theta_{{\alpha_2},{\alpha_4}}$ would contain $x_{\alpha_1}^l$ for some $1\le l< k$, but this is a contradiction since any element of $I$ cannot contain $x_{\alpha_1}^l$. In the remaining case $(0,0,1)$, we have $(a_1,a_2)=(\frac{1}{2},-\frac{1}{2})$. If $\{x_{\alpha_3},\Theta_{{\alpha_2},{\alpha_4}}\}$ contains $x_{\alpha_1}^{k-1}x_{\alpha_4}$, we similarly have $(a_1,a_2)=(-\frac{1}{2},\frac{1}{2})$. 

We can use the same argument for $J_{\alpha_2, \alpha_3, {\alpha_5}}$, and we can conclude that there are two possibilities as follows:\\
$\bullet$ $(a_1,a_2)=(\frac{1}{2},-\frac{1}{2})$; $\Theta_{{\alpha_3},{\alpha_4}}$ contains $x_{\alpha_1}^l x_{\alpha_5}$, and $\Theta_{{\alpha_2},{\alpha_5}}$ contains $x_{\alpha_1}^m x_{\alpha_4}$ with $l,m\ge0$ with $l+m=k-1$.\\
$\bullet$ $(a_1,a_2)=(-\frac{1}{2},\frac{1}{2})$; $\Theta_{{\alpha_2},{\alpha_4}}$ contains $x_{\alpha_1}^l x_{\alpha_5}$, and $\Theta_{{\alpha_3},{\alpha_5}}$ contains $x_{\alpha_1}^m x_{\alpha_4}$  with $l,m\ge0$ with $l+m=k-1$.\\

Using the condition $w(\Theta_{i,j})=w_i+w_j$ in Lemma \ref{wt}, we see that $w(\Theta_{i,j})\in\frac{1}{2}+\mathbb{Z}$ for $i=\alpha_1,\alpha_2,\alpha_3$ and $j=\alpha_4,\alpha_5$. This implies that these $\Theta_{i,j}$'s are in the ideal $(x_{\alpha_4},x_{\alpha_5})$ since $w(x_{\alpha_i})\in\mathbb{Z}$ for $i=1,2,3$. Similarly, every term of $\Theta_{{\alpha_1},{\alpha_2}},\Theta_{{\alpha_1},{\alpha_3}},\Theta_{{\alpha_2},{\alpha_3}}$ and $\Theta_{{\alpha_4},{\alpha_5}}$ is either in $\mathbb{C}[x_{\alpha_1},x_{\alpha_2},x_{\alpha_3}]$ or in $(x_{\alpha_4},x_{\alpha_5})^2$. By (\ref{4-dim}), we see that $p(x_1,x_2,x_3,x_4,x_5)$ in (\ref{f}) is contained in $(x_{\alpha_4},x_{\alpha_5})^2$. Then we can write
$$\{x_{\alpha_i},f\}=\{x_{\alpha_i},h(x_{\alpha_1},x_{\alpha_2},x_{\alpha_3})\}\Delta^2+p_{(i)}(x_1,x_2,x_3,x_4,x_5)$$
with $p_{(i)}\in(x_{\alpha_4},x_{\alpha_5})$ for $i=1,2,3$. Note that $\{x_{\alpha_i},\Delta\}$ is in $(x_{\alpha_4},x_{\alpha_5})^2$ for $i=1,2,3$. We have $\{x_{\alpha_i},f\}=0$ by (\ref{Pf}), and therefore $\{x_{\alpha_i},h\}$ is in $(x_{\alpha_4},x_{\alpha_5})$ for $i=1,2,3$. By the condition $\{x_{\alpha_1},h\}\in(x_{\alpha_4},x_{\alpha_5})$, we have $h\in \mathbb{C}[x_{\alpha_1},\Delta]$. We can also show that $h=\kappa\Delta^{n-2}$ with $n\ge2$ for some nonzero constant $\kappa$ by using $\{x_{\alpha_2},h\}\in(x_{\alpha_4},x_{\alpha_5})$.


We show that $l=m=0$. First we consider the case when $(a_1,a_2)=(\frac{1}{2},-\frac{1}{2})$. Then we see that $\Theta_{{\alpha_3},{\alpha_4}}$ contains $x_{\alpha_1}^l x_{\alpha_5}$ and other terms of $\Theta_{{\alpha_3},{\alpha_4}}$ are in $(x_{\alpha_4},x_{\alpha_5})^3$ since we have $\deg(\Theta_{{\alpha_3},{\alpha_4}})<\deg(x_{\alpha_3}x_{\alpha_4})$ and $w(\Theta_{{\alpha_3},{\alpha_4}})=-\frac{1}{2}$. (Note that $\Theta_{{\alpha_3},{\alpha_4}}$ does not contain $x_{\alpha_2}x_{\alpha_3} x_{\alpha_5}$ since $\deg(x_{\alpha_1}^l)<\deg(x_{\alpha_1}^{k+1})=\deg(x_{\alpha_2}x_{\alpha_3})$.) This implies that $\{x_{\alpha_2},\Theta_{{\alpha_3},{\alpha_4}}\}$ contains $x_{\alpha_1}^{l-1}x_{\alpha_2}x_{\alpha_5}$ unless $l=0$. We show that this is contrary to $J_{{\alpha_2},{\alpha_3},{\alpha_4}}=0$. It suffices to show that $\{x_{\alpha_3},\Theta_{{\alpha_2},{\alpha_4}}\}$ and $\{x_{\alpha_4},\Theta_{{\alpha_2},{\alpha_3}}\}$ do not contain $x_{\alpha_1}^{l-1}x_{\alpha_2}x_{\alpha_5}$. In order for $\{x_{\alpha_3},\Theta_{{\alpha_2},{\alpha_4}}\}$ not to contain $x_{\alpha_1}^{l-1}x_{\alpha_2}x_{\alpha_5}$, it suffices to show that $\Theta_{{\alpha_2},{\alpha_4}}$ is in $(x_{\alpha_4},x_{\alpha_5})^2$ since $\{x_3,\Theta_{{\alpha_2},{\alpha_4}}\}$ would be in $(x_{\alpha_4},x_{\alpha_5})^2$ in this case. By considering the two degrees of $\Theta_{{\alpha_2},{\alpha_4}}$, the only possible monomial of $\Theta_{{\alpha_2},{\alpha_4}}$ which is not in $(x_{\alpha_4},x_{\alpha_5})^2$ is $x_{\alpha_2}^2x_{\alpha_5}$. However, the condition that $\Theta_{{\alpha_2},{\alpha_5}}$ contains $x_{\alpha_1}^m x_{\alpha_4}$ implies that this is impossible since
$$\deg(x_{\alpha_2}^2x_{\alpha_5})>2d_{\alpha_2}+d_{\alpha_5}-s=d_{\alpha_2}+md_{\alpha_1}+d_{\alpha_4}>\deg(\Theta_{{\alpha_2},{\alpha_4}}).$$
As for $\{x_{\alpha_4},\Theta_{{\alpha_2},{\alpha_3}}\}$, first note that $\Theta_{{\alpha_2},{\alpha_3}}-x_{\alpha_1}^k$ is in $(x_{\alpha_4},x_{\alpha_5})^2$. We see that, in order for $\{x_{\alpha_4},\Theta_{{\alpha_2},{\alpha_3}}\}$ to contain $x_{\alpha_1}^{l-1}x_{\alpha_2}x_{\alpha_5}$, the polynomial $\Theta_{{\alpha_2},{\alpha_3}}-x_{\alpha_1}^k$ must contain $x_{\alpha_1}^{l'}x_{\alpha_2}x_{\alpha_5}^2$ for some $0\le l'\le l-1$ and $\Theta_{{\alpha_4},{\alpha_5}}$ must contain $x_{\alpha_1}^{l-l'-1}$ by considering the two degrees of $\Theta_{{\alpha_2},{\alpha_3}}$. We have $l+m=k-1$ and in particular $l-l'-1\le k-1$. However, (\ref{f}) shows that $\deg(f)\ge 2\deg(\Delta)=2(k+1)d_{\alpha_1}$, which implies that
$$d_{\alpha_4}+d_{\alpha_5}\ge (k+2)d_{\alpha_1}.$$
Thus, $\Theta_{{\alpha_4},{\alpha_5}}$ cannot contain $x_{\alpha_1}^{l-l'-1}$. Therefore, we have shown that $\{x_{\alpha_4},\Theta_{{\alpha_2},{\alpha_3}}\}$ does not contain $x_{\alpha_1}^{l-1}x_{\alpha_2}x_{\alpha_5}$, and we obtain $l=0$. Similarly we can show that $m=0$ using $J_{{\alpha_2},{\alpha_3},{\alpha_5}}=0$. The same argument works for the case when $(a_1,a_2)=(-\frac{1}{2},\frac{1}{2})$ by switching the roles of $x_{\alpha_4}$ and $x_{\alpha_5}$. Note that we did not use $d_{\alpha_4}\le d_{\alpha_5}$ in the above argument.

When $(a_1,a_2)=(\frac{1}{2},-\frac{1}{2})$, we have $k=l+m+1=1$. Since $\deg(f)=\deg(\Delta^n)=2ns$ and $\deg(f)=\sum_{i=1}^5 d_i-2s$ (see (\ref{4-dim})), the degrees $d_i$ are written as
\begin{equation}\label{degree}
d_{\alpha_2}=s-2t,\;d_{\alpha_1}=s,\;d_{\alpha_3}=s+2t,\;d_{\alpha_4}=\left(n-\frac{1}{2}\right)s-t,\;d_{\alpha_5}=\left(n-\frac{1}{2}\right)s+t
\end{equation}
for some $0\le t<\frac{1}{2}s$. Now we can deduce the Poisson matrix $\Theta_n$ in Section \ref{2}. By considering possible degrees and the condition $w(\Theta_{i,j})=w_i+w_j$, we can determine the entries of the Poisson matrix except $\Theta_{{\alpha_4},{\alpha_5}}$. We see that $\Theta_{{\alpha_4},{\alpha_5}}\in\mathbb{C}[x_{\alpha_1},x_{\alpha_2},x_{\alpha_3}]$ for degree reasons and can show that it is equal to $n\kappa\Delta^{n-1}$ from (\ref{4-dim}) and $h=\kappa\Delta^{n-1}$. Thus, we finally obtain the following Poisson matrix
$$\begin{pmatrix}
0&x_{\alpha_2}&-x_{\alpha_3}&\frac{1}{2}x_{\alpha_4}&-\frac{1}{2}x_{\alpha_5}\\
-x_{\alpha_2}&0&x_{\alpha_1}&0&c_1 x_{\alpha_4}\\
x_{\alpha_3}&-x_{\alpha_1}&0&c_2 x_{\alpha_5}&0\\
-\frac{1}{2}x_{\alpha_4}&0&-c_2 x_{\alpha_4}&0&n\kappa \Delta^{n-1}\\
\frac{1}{2}x_{\alpha_5}&-c_1 x_{\alpha_4}&0&-n\kappa \Delta^{n-1}&0
\end{pmatrix}$$
with nonzero constants $c_1,c_2$ and $\kappa$ where the ordering of the variables for this matrix is given by $x_{\alpha_1},x_{\alpha_2},x_{\alpha_3},x_{\alpha_4},x_{\alpha_5}$. From $J_{{\alpha_2},{\alpha_3},{\alpha_4}}=0$, we have $c_1 c_2=\frac{1}{2}$.

We obtain the matrix $\Theta_n$ in Section \ref{2} by replacing $x_{\alpha_1},x_{\alpha_2},x_{\alpha_3},x_{\alpha_4}$ and $x_{\alpha_5}$ with $\frac{1}{2}h,\frac{1}{2}x,y,\frac{\sqrt{\kappa}}{2}e_0$ and $\frac{\sqrt{\kappa}}{c_1}e_1$ respectively. For example, we have
$$\{x,h\}=\{2x_{\alpha_2},2x_{\alpha_1}\}=-4x_{\alpha_2}=-2x.$$
When $(a_1,a_2)=(-\frac{1}{2},\frac{1}{2})$, we have $d_{\alpha_4}=d_{\alpha_5}$ and we are reduced to the former case by switching the roles of $x_{\alpha_4}$ and $x_{\alpha_5}$ again. Note that we have $t=0$ in this case since $d_{\alpha_4}=d_{\alpha_5}$.

\vspace{3mm}

\noindent
Case (II): $a_1=-a_2=0$\\
We show that this case does not happen.

By the indecomposability of $X$, the radical of the ideal $\mathfrak{a}$ generated by $\Theta_{i,j}$'s must coincide with the maximal ideal $(x_1,\dots,x_5)$. In particular, $\mathfrak{a}$ and hence some $\Theta_{i,j}$ must contain powers of $x_{\alpha_5}$. Since $w_{\alpha_5}=0$, a power of $x_{\alpha_5}$ can only locate at $\Theta_{\alpha_1,\alpha_4},\;\Theta_{\alpha_1,\alpha_5},\;\Theta_{{\alpha_2},{\alpha_3}}$ or $\Theta_{{\alpha_4},{\alpha_5}}$. Note that $x_5$ cannot locate at $\Theta_{\alpha_1,\alpha_4}$ or $\Theta_{\alpha_1,\alpha_5}$ by the choice of the coordinates and by the assumption that $a_2=0$. Squares or higher powers of $x_5$ also cannot locate at $\Theta_{\alpha_1,\alpha_4}$ or $\Theta_{\alpha_1,\alpha_5}$ since
$$\deg(\Theta_{\alpha_1,\alpha_4})=d_{\alpha_4}<d_{\alpha_5}=\deg(\Theta_{\alpha_1,\alpha_5})<2d_{\alpha_5}.$$
If $\Theta_{{\alpha_4},{\alpha_5}}$ contains a power $x_{\alpha_5}^l$, then $l$ is 1 since $\alpha_4<\alpha_5$. Also, $d_{\alpha_4}=s$ and we have $\Theta_{{\alpha_i},{\alpha_4}}=0$ for $i=1,2,3$ since $\deg(\Theta_{{\alpha_i},{\alpha_4}})=d_{\alpha_i},\,w(\Theta_{{\alpha_i},{\alpha_4}})=w_{\alpha_i}$ and $\Theta_{{\alpha_i},{\alpha_4}}\in I$. In this case $\frac{\partial f}{\partial x_{\alpha_5}}=0$ by (\ref{4-dim}), which is a contradiction \cite[Lemma 2.6(2)]{LNSvS}.

Therefore, a power $x_{\alpha_5}^l,\,l\ge1$ locates at $\Theta_{{\alpha_2},{\alpha_3}}$. Then we have $d_{\alpha_2}+d_{\alpha_3}-s=l d_{\alpha_5}$, and thus
$$d_{\alpha_4}\le d_{\alpha_5}=\frac{1}{l}(d_{\alpha_2}+d_{\alpha_3}-s)< d_{\alpha_2}+d_{\alpha_3}=\deg(\Delta).$$
Then we have $\Theta_{{\alpha_i},{\alpha_4}}=0$ for $i=1,2,3$. Indeed, we have $w(\Theta_{{\alpha_i},{\alpha_4}})=w_{\alpha_i}$ and $\Theta_{{\alpha_i},{\alpha_4}}\in I$, but the degree of $\Theta_{{\alpha_i},{\alpha_4}}$ is less than $\deg(x_{\alpha_i}x_{\alpha_4})$ and $\deg(x_{\alpha_i}\Delta)$, which forces $\Theta_{{\alpha_i},{\alpha_4}}$ to be zero. Therefore, we have the same contradiction as above.
\qed

\begin{center}
Kavli Institute for the Physics and Mathematics of the Universe (WPI), UTIAS,\\
The University of Tokyo, Kashiwa, Chiba, 277-8583, Japan.\\
e-mail address: ryo.yamagishi@ipmu.jp
\end{center}

\section*{Appendix: Klein singularities from contact point of view, by Yoshinori Namikawa}

Let $R = \oplus_{i \geq 0}R_i$ be a positively graded, finitely generated $\mathbf{C}$-algebra.  
We assume that $R_0 = \mathbf{C}$ and $R$ is an integrally closed domain. Then the affine 
normal variety $X := \mathrm{Spec} \; R$ has a $\mathbf{C}^*$-action with a unique 
fixed point at the origin $0 \in X$. When $X$ admits an algebraic symplectic 2-form $\omega$
on the smooth part $X_{reg}$ so that $\omega$ is homogeneous with respect to the 
$\mathbf{C}^*$-action and $\omega$ extends to a regular 2-form on a resolution 
$f: Y \to X$ of $X$, we call the pair $(X, \omega)$ a conical symplectic variety.  

The aim of this short note is to explain the following result from 
contact point of view. In the below $\mathbf{C}^2(x,y)$ means the 2-dimensional affine 
space with coordinates $x$ and $y$.  

{\bf Proposition}. {\em Let $(X, \omega)$ be a conical symplectic variety 
of dimension $2$. Then $(X, \omega)$ is isomorphic to       
one of the following} 
\vspace{0.2cm}

(i) Smooth case: $X = \mathbf{C}^2(x, y), \;\ \omega = dx \wedge dy$ $wt(x, y) = (a, b)$, where $a$ and $b$ are relatively prime 
positive integers.    
 
(ii) $A_{n-1}$-type $(n \geq 2)$: $X = \{xy - z^n = 0\} \subset \mathbf{C}^3(x, y, z), \; \; 
\omega = Res (dx \wedge dy \wedge dz/xy - z^n)$, $wt(x,y,z) = (a, b, c)$ where $a$, $b$ and $c$ are positive integers such that $\mathrm{g.c.m.}(a, b, c) = 1$ and  
$a + b = nc $.  

(iii) $D_n$-type $(n \geq 4)$:  $X = \{x^{n-1} + xy^2 + z^2 = 0\} \subset \mathbf{C}^3(x,y,z), \;\; 
\omega = Res (dx \wedge dy \wedge dz/x^{n-1} + xy^2 + z^2)$, $wt (x,y,z) = (2, n-2, n-1)$. 

(iv) $E_6$-type: $X = \{x^4 + y^3 + z^2 = 0\} \subset \mathbf{C}^3(x,y,z), \;\; 
\omega = Res (dx \wedge dy \wedge dz/x^4 + y^3 + z^2)$, $wt (x,y,z) = (3, 4, 6)$   
       
(v) $E_7$-type: $X = \{x^3y + y^3 + z^2 = 0\} \subset \mathbf{C}^3(x,y,z), \;\; 
\omega = Res (dx \wedge dy \wedge dz/x^3y + y^3 + z^2)$, $wt (x,y,z) = (4, 6, 9)$

(vi) $E_8$-type: $X = \{x^5 + y^3 + z^2 = 0\} \subset \mathbf{C}^3(x,y,z), \;\; 
\omega = Res (dx \wedge dy \wedge dz/x^5 + y^3 + z^2)$, $wt (x,y,z) = (6, 10, 15)$ 
\vspace{0.2cm}

The relationship between a conical symplectic variety and a contact Fano orbifold is 
the following (for details, see [Na 1, 4.4], [Na 2] and [Na 3, \S. 2]). Let  
$(X, \omega)$ be a conical symplectic variety of dimension $2d$ and put 
$l := \mathrm{deg}(\omega)$. We know that $l > 0$. Then we take a minimal homogeneous generator $x_0$, ..., $x_n$ of the $\mathbf{C}$-algebra $R$ and put $a_i := \mathrm{deg} (x_i)$. 
We may assume that $\{a_0, ..., a_n\}$ have no common factors. In fact, if they have a common factor, say $f$, we can put a new degree for each element of $x \in R$ by $\mathrm{deg}(x)/f$. 
$X$ is embedded in an affine space $\mathbf{C}^{n+1}$ by using $x_i$'s. The affine space 
$\mathbf{C}^{n+1}$ has a $\mathbf{C}^*$-action determined by $(x_0, ..., x_n) \to 
(t^{a_0}x_0, ..., t^{a_n}x_n)$ for $t \in \mathbf{C}^*$. The quotient space $\mathbf{C}^{n+1} - \{0\}/\mathbf{C}^*$ is nothing but the weighted projective space $\mathbf{P}(a_0, ..., a_n)$. 
Here we put $\mathbf{P}(X) := X - \{0\}/\mathbf{C}^*$. Then $\mathbf{P}(X)$ is a subvariety 
of $\mathbf{P}(a_0, ..., a_n)$. Notice that $\mathbf{P}(a_0, ..., a_n)$ has a natural orbifold 
structure (hereafter we denote such an orbifold by $\mathbf{P}(a_0, ..., a_n)^{orb}$) and a tautological orbifold line bundle $\mathcal{O}_{\mathbf{P}(a_0, ..., a_n)^{orb}}(1)$. 
They induces an orbifold structure $\mathbf{P}(X)^{orb}$ on $\mathbf{P}(X)$ and a tautological orbifold line bundle 
$\mathcal{L}$ on $\mathbf{P}(X)^{orb}$. Then the orbifold 
$\mathbf{P}(X)^{orb}$ has a contact structure with a contact orbifold line bundle $$\mathcal{M} := \mathcal{L}^{\otimes l}.$$ As in the usual situation we have an isomorphism $$-K_{\mathbf{P}(X)^{orb}} \cong \mathcal{M}^{\otimes d}.$$   
In order to recover the original conical symplectic variety $(X, \omega)$, we just put 
$$X := \mathrm{Spec}\: \oplus_{i \geq 0} H^0(\mathbf{P}(X)^{orb}, \mathcal{L}^{\otimes i}).$$  
We can define a homogeneous symplectic structure $\omega$ by using the orbifold contact structure on $\mathbf{P}(X)^{orb}$. 

The orbifold structure on $\mathbf{P}(X)$ determines a finite number of prime Weil divisor $\{D_{\alpha}\}$; each of them is called a {\em ramification divisor} and it is attached with a positive integer $e_{\alpha} > 1$ called the 
ramification index. Now one can consider a $\mathbf{Q}$-divisor 
$\Delta$ on $\mathbf{P}(X)$ by 
$$\Delta := \sum_{\alpha}  (1 - \frac{1}{e_{\alpha}})D_{\alpha}.$$ By Lemma 2.1 of [Na 2], 
the log pair $(\mathbf{P}(X), \Delta)$ is a log Fano variety. In particular, 
$-(K_{\mathbf{P}(X)} + \Delta)$ is an ample $\mathbf{Q}$-divisor.  
 
When $d = 1$ (that is, when $\mathbf{P}(X)$ is a curve), the situation is extremely simple.  
Since $-(K_{\mathbf{P}(X)} + \Delta)$ is an ample $\mathbf{Q}$-divisor, one has  
$\mathbf{P}(X) \cong \mathbf{P}^1$. Then, since $\mathrm{deg}\; K_{\mathbf{P}^1} = -2$, the inequality 
$$ -2 + \sum_{\alpha} (1 - \frac{1}{e_{\alpha}}) < 0 $$ must hold. 
Since each term $1 - \frac{1}{e_{\alpha}} \geq \frac{1}{2}$, the number of the ramification 
points does not exceed $3$. The following are possible candidates of the ramification 
indexes: 
$$e_1 = a \;\; a \geq 2 $$
$$ (e_1, e_2) = (a, b) \;\; a \geq 2, \;\;  b \geq 2$$
$$ (e_1, e_2, e_3) = (2, 2, n), \;\; n \geq 2$$
$$ (e_1, e_2, e_3) = (2, 3, 3), \;\;  (2,3,4), \;\; (2,3,5) $$     
As any different 3 points on $\mathbf{P}^1$ can be transformed to 
$\{0, \;\;1,\;\;  \infty\}$, we see that these numerical data completely 
determine the orbifold structure. We recommend the reader to check the following.  

{\bf Lemma}.  {\em In the list of Proposition, the cases (i) and (ii) give the orbifolds with 
$(e_1, e_2) = (a, b)$. Here, when $(e_1, e_2) = (a, 1)$ or $(1, b)$, we interpret it respectively as 
$e_1 = a$ or $e_1 = b$. The case (iii) gives the orbifold with $(e_1, e_2, e_3) = (2,2, n-2)$. 
The case (iv) gives the orbifold with $(e_1, e_2, e_3) = (2,3,3)$. The case (v) gives 
the orbifold with $(e_1, e_2, e_3) = (2,3,4)$, and finally, the case (vi) gives the orbifold 
with $(e_1, e_2, e_3) = (2, 3, 5)$.}   \vspace{0.2cm}    

Notice here that in the cases (i) and (ii), different conical symplectic varieties give 
the same orbifold structure. In order to specify a conical symplectic variety, we must additionally fix a tautological orbifold line bundle on the orbifold. Here we introduce two invariants of an orbifold line bundle. Let $\mathbf{P}^{1, orb}$ be an orbifold whose underlying space is $\mathbf{P}^1$. Let $\mathcal{L}$ be an 
orbifold line bundle on $\mathbf{P}^{1, orb}$. 
Let $Q_1$, ..., $Q_r$ $(r \leq 3)$ be ramification points of $\mathbf{P}^{1, orb}$ with ramification indexes $e_1$, ..., $e_r$.  Around each $Q_i$, $\mathcal{L}$ is not necessarily 
a usual line bundle. Now let $\tau_i({\mathcal L})$ be the smallest positive integer $k$ 
such that $\mathcal{L}^{\otimes k}$ is a usual line bundle around $Q_i$. By definition 
$\tau_i({\mathcal L})$ is a divisor of $e_i$ and 
$$\tau_i(\mathcal{L}^{\otimes m}) = \frac{\tau_i(\mathcal{L})}{\mathrm{g.c.d.}(m, \tau_i(\mathcal{L}))}$$ 
for any positive integer $m$. An orbifold line bundle 
$\mathcal{L}$ is a usual line bundle if and only if $\tau_i(\mathcal{L}) = 1$ for all 
$i$. Let $e$ be the least common multiple of $e_1$, ... $e_r$. Then $\mathcal{L}^{\otimes e}$ 
is a usual line bundle on $\mathbf{P}^1$. We define 
$$ \mathrm{deg}\; \mathcal{L} := \frac {\mathrm{deg}\; \mathcal{L}^{\otimes e}}{e}.$$
By definition $$\mathrm{deg}\; \mathcal{L} \in \frac{1}{e} \mathbf{Z}$$ for any orbifold 
line bundle $\mathcal{L}$. 
For example, we have 
$$ \tau_i (K_{\mathbf{P}^{1, orb}}) = e_i, \; \; (\mathrm{hence}, \;\; \tau_i (-K_{\mathbf{P}^{1, orb}}) = e_i)$$
$$\mathrm{deg}\; K_{\mathbf{P}^{1, orb}} = -2 + \sum_{1 \le i \le r} (1 - \frac{1}{e_i}).$$  
Let us consider how many different tautological orbifold line bundles $\mathbf{P}^{1, orb}$ has. 
When $(e_1, e_2, e_3) = (2, 3, 5)$, $\mathrm{deg}\; K_{\mathbf{P}^{1,orb}} = -\frac{1}{30}$. 
The degree of any orbifold line bundle is in $\frac{1}{30}\mathbf{Z}$ because $\mathrm{l.c.m}(2,3,5) = 30$. This means that the tautological line bundle must be $-K_{\mathbf{P}^{1,orb}}$ itself.  In the similar ways, one can check that the tautological line bundle must be $-K_{\mathbf{P}^{1,orb}}$ when $(e_1, e_2, e_3) = (2,3,3)$, $(2,3, 4)$ and $(2, 2, n)$ with 
$n$ even. When $(e_1, e_2, e_3) = (2, 2, n)$ with $n$ odd, we have $\mathrm{l.c.m}(2,2,n) = 2n$; hence the degree of an arbitrary orbifold line bundle is an element of $\frac{1}{2n}\mathbf{Z}$.  On the other hand, $\mathrm{deg}\; K_{\mathbf{P}^{1,orb}} = -\frac{1}{n}$. We show that 
there is no orbifold line bundle $\mathcal{L}$ with $\mathcal{L}^{\otimes 2} \cong   
-K_{\mathbf{P}^{1,orb}}$. In fact, we always have $\tau_1(\mathcal{L}^{\otimes 2}) = 1$ for any 
orbifold line bundle $\mathcal{L}$, but $\tau_1(-K_{\mathbf{P}^{1, orb}}) = 2$.

When $(e_1, e_2) = (a, b)$, we take the greatest common divisor $m$ of $a$ and $b$, and 
write $a = ma'$, $b = mb'$. In this case $\mathrm{deg}\; K_{\mathbf{P}^{1,orb}} = -\frac{a' + b'}
{ma'b'}$. Let $c$ be a divisor of $a' + b'$ and suppose that an orbifold line bundle $\mathcal{L}$ 
satisfies $\mathcal{L}^{\otimes c} \cong -K_{\mathbf{P}^{1, orb}}$. Since $\tau_1(-K_{\mathbf{P}^{1, orb}}) = a$ and $\tau_2(-K_{\mathbf{P}^{1, orb}}) = b$, we 
see that $\mathrm{g.c.d.}(c, a) = \mathrm{g.c.d.}(c, b) = 1$. 
We first prove that the following two conditions are the same. 

(i) $c$ is a divisor of $a' + b'$ with $\mathrm{g.c.d.}(c, a) = \mathrm{g.c.d.}(c, b) = 1$. 

(ii) $c$ is a divisor of $a + b$ with $\mathrm{g.c.d.}(c, a, b) = 1$. 

Since it is clear that (i) implies that (ii), we show that (ii) implies that (i). 
Take $c$ in (ii). Then  $\mathrm{g.c.d.}(c, a) = \mathrm{g.c.d.}(c, b) = 1$. 
In fact, if $d$ is a common divisor of $c$ and $a$ (resp. $c$ and $b$), then $d$ is also a divisor of $b$ (resp. $a$) because $c$ is a divisor of $a + b$.
On the other hand, $c$ must be a divisor of $a' + b'$. In fact, if $c$ is not a divisor of $a' + b'$, then $\mathrm{g.c.d.}(c, m) \ne 1$. But 
this contradicts $\mathrm{g.c.d.}(a, b, c) = 1$. Therefore (ii) implies (i). 

We next prove that, for any $c$ in (i) (or (ii)), the orbifold line bundle $\mathcal{L}$ with $\mathcal{L}^{\otimes c} \cong -K_{\mathbf{P}^{1, orb}}$ must be unique if it exists. 
If another orbifold line bundle $\mathcal{L}'$ satisfies this property, then 
$(\mathcal{L} \otimes (\mathcal{L}')^{-1})^{\otimes c} \cong \mathcal{O}_{\mathbf{P}^1}$. 
We write $M = \mathcal{L} \otimes (\mathcal{L}')^{-1}$ for simplicity. Then 
$\tau_1(M^{\otimes c}) = \tau_2(M^{\otimes c}) = 1$. 
Since $\mathrm{g.c.d.}(a, c) = \mathrm{g.c.d.}(b, c) = 1$, we have 
$\tau_1(M) = \tau_2 (M) = 1$. This implies that $M$ is a usual line bundle on $\mathbf{P}^1$. 
On the other hand, $\mathrm{deg}\; M = 0$. Hence $M \cong \mathcal{O}_{\mathbf{P}^1}$. 
 
Finally we show that, for any $c$ in (ii), there exists an orbifold line bundle  $\mathcal{L}$ on $\mathbf{P}^{1, orb}$ (= tautological line 
bundle) such that $\mathcal{L}^{\otimes c} \cong 
-K_{\mathbf{P}^{1, orb}}$ and $\mathbf{P}^{1, orb}$ has a contact orbifold structure with contact line bundle $-K_{\mathbf{P}^{1, orb}}$. To construct such an orbifold structure, we just 
take an $X$ with $a \geq 2$ and $b \geq 2$ in the list (ii) of Proposition, and consider $\mathbf{P}(X)^{orb}$. The $c$ in the argument above is nothing but the $c$ in the list (ii) of Proposition.  

When $e_1 = a$, we have $\mathrm{deg}\; K_{\mathbf{P}^{1,orb}} = - \frac{a + 1}{a}$. 
Notice that $\mathrm{deg}\; \mathcal{L} \in \frac{1}{a} \mathbf{Z}$ for any orbifold line bundle 
$\mathcal{L}$. Let $c$ be a divisor $c$ of $a + 1$ and consider the existence problem of 
an orbifold $\mathcal{L}$ with  
$\mathcal{L}^{\otimes c} \cong -K_{\mathbf{P}^{1, orb}}$.  For each $c$, if such an $\mathcal{L}$ exists, then it is unique by the same argument in the case 
$(e_1, e_2) = (a, b)$. Finally the existence of such an $\mathcal{L}$ (for any divisor $c$ of $a + 1$) follows by taking $\mathbf{P}(X)^{orb}$ for an $X$ with $a = 1$ in the list (ii) of Proposition or an $X$ in the list (i) of Proposition. \vspace{0.2cm}  

\begin{center}
{\bf References} 
\end{center}
\vspace{0.2cm}

[Na 1] Namikawa, Y.: Equivalence of symplectic singularities, Kyoto Journal of Math. {\bf 53} 
(2013) 483-514 

[Na 2] Namikawa, Y.: Fundamental groups of symplectic singularities, Adv. Stud. Pure Math. 
{\bf 74} (2017) Higher dimensional algebraic geometry 321-334

[Na 3] Namikawa, Y.: A finiteness theorem on symplectic singularities, Compos. Math. {\bf 
152} (2016) 1225-1236  

\vspace{3mm}

\begin{center}
Department of Mathematics, Faculty of Science, Kyoto University, Japan\\
namikawa@math.kyoto-u.ac.jp
\end{center}


\begin{thebibliography}{99}
\bibitem[B]{B}
A. Beauville, Symplectic singularities. Invent. Math. 139 (2000), no. 3, 541--549.

\bibitem[K1]{K1}
D. Kaledin, Symplectic singularities from the Poisson point of view, J. Reine Angew. Math. {\bf 600} (2006), 135--156. 

\bibitem[K2]{K2}
D. Kaledin, Geometry and topology of symplectic resolutions, in {\it Algebraic geometry---Seattle 2005. Part 2}, 595--628, Proc. Sympos. Pure Math., 80, Part 2, Amer. Math. Soc., Providence, RI.

\bibitem[K3]{K3}
D. B. Kaledin, Normalization of a Poisson algebra is Poisson. Proc. Steklov Inst. Math. {\bf 264} (2009), no.~1, 70--73; translated from Tr. Mat. Inst. Steklova {\bf 264} (2009), Mnogomernaya Algebraicheskaya Geometriya, 77--80. 

\bibitem[LNS]{LNS}
M. Lehn, Y. Namikawa, C. Sorger: Slodowy slices and universal Poisson deformations.
Compositio Math. 148 (2012), 121 - 144.

\bibitem[LNSvS]{LNSvS}
M. Lehn, Y. Namikawa, C. Sorger, D. van Straten, On symplectic hypersurfaces, in {\it Minimal models and extremal rays (Kyoto, 2011)}, 277--298, Adv. Stud. Pure Math., 70, Math. Soc. Japan.

\bibitem[N1]{N1}
Y. Namikawa, On the structure of homogeneous symplectic varieties of complete intersection, Invent. Math. {\bf 193} (2013), no.~1, 159--185.

\bibitem[N2]{N2}
Y. Namikawa, A finiteness theorem on symplectic singularities, Compos. Math. {\bf 152} (2016), no.~6, 1225--1236. 

\bibitem[P]{P}
M. H. Park, Integral closure of graded integral domains, Comm. Algebra {\bf 35} (2007), no.~12, 3965--3978.
\end{thebibliography}
\end{document}